\documentclass[11pt]{amsart}
\usepackage{amsmath}
\usepackage{amssymb}
\usepackage{amsfonts}
\usepackage{amsthm}
\usepackage{enumerate}
\usepackage[mathscr]{eucal}
\usepackage{verbatim}
\usepackage{amsthm}
\usepackage{amscd}
\usepackage[mathscr]{eucal}
\usepackage{appendix}
\usepackage{tikz}
\usepackage{hyperref}
\usepackage{cleveref}
\usepackage[normalem]{ulem}

\usepackage{xcolor}
\hypersetup{
    colorlinks,
    linkcolor={blue!50!black},
    citecolor={blue!50!black},
    urlcolor={blue!80!black}
}

\numberwithin{equation}{section}

\newtheorem{innercustomgeneric}{\customgenericname}
\providecommand{\customgenericname}{}

\newcommand{\newcustomtheorem}[2]{\newenvironment{#1}[1]
  {\renewcommand\customgenericname{#2}
   \renewcommand\theinnercustomgeneric{##1}\innercustomgeneric}{\endinnercustomgeneric}}

\newcustomtheorem{customthm}{Theorem}

\newcommand{\newcustomlemma}[2]{\newenvironment{#1}[1]
  {\renewcommand\customgenericname{#2}
   \renewcommand\theinnercustomgeneric{##1} \innercustomgeneric}{\endinnercustomgeneric}}

\newcustomlemma{customlemma}{Lemma}

\newcustomlemma{customproposition}{Proposition}

\newcustomlemma{customclaim}{Claim}

\theoremstyle{plain}
\newtheorem{theorem}{Theorem}
\newtheorem{corollary}[theorem]{Corollary}
\newtheorem{lemma}[theorem]{Lemma}
\newtheorem{proposition}[theorem]{Proposition}
\newtheorem{remark}{Remark}

\newtheorem*{theorem*}{Theorem}
\newtheorem*{lemma*}{Lemma}
\newtheorem*{proposition*}{Proposition}
\newtheorem*{corollary*}{Corollary}
\newtheorem*{remark*}{Remark} 
\newtheorem*{remarks*}{Remarks}
\newtheorem*{conj*}{Conjecture}

\def\N{{\mathbb N}}

\def\S{\mathbb{S}}

\def\R{{\mathbb R}}
\def\Z{{\mathbb Z}}

\newcommand{\bbz}{\mathbb{Z}}

\newcommand{\bbs}{\mathbb S}

\newcommand{\bbr}{\mathbb{R}}
\newcommand{\bbrn}{\mathbb R^n}

\newcommand{\bbn}{\mathbb{N}}

\newcommand{\MM}{\mathcal{M}}

\newcommand{\xxxi}{\vec{\boldsymbol{\xi}\;}}
\newcommand{\xxx}{\vec{\boldsymbol{x}}}
\newcommand{\yyy}{\vec{\boldsymbol{y}}}

\newcommand{\fff}{\vec{\boldsymbol{f}}}

\newcommand{\xxi}{\vec{\boldsymbol{\xi}}}

\def\000{\vec{\boldsymbol{0}}}

\newcommand{\q}{\quad}
\newcommand{\qq}{\qquad}

\DeclareFontFamily{U}{mathx}{\hyphenchar\font45}
\DeclareFontShape{U}{mathx}{m}{n}{
	<5> <6> <7> <8> <9> <10>
	<10.95> <12> <14.4> <17.28> <20.74> <24.88>
	mathx10
}{}

\newcommand{\widecheck}{^\vee}

\def\wh{\widehat}

\newcommand{\todo}[1]{{\color{red}{[#1]}}}

\newcommand{\wt}{\widetilde}

\newcommand{\supp}{\mathrm{supp}}

\makeatletter
\@namedef{subjclassname@2020}{\textup{2020} Mathematics Subject Classification}
\makeatother

\makeindex

\begin{document}

		\author{Georgios~Dosidis}
\address{Department of Mathematical Analysis, Faculty of Mathematics and Physics, Charles University, Sokolovsk\'a 83, 186 75 Praha 8, Czech Republic}
\email{Dosidis@karlin.mff.cuni.cz}
		
		\author{Bae Jun Park}
\address{Department of Mathematics, Sungkyunkwan University, Suwon 16419, Republic of Korea}
\email{bpark43@skku.edu}
		
		\author{Lenka~Slav{\'i}kov{\'a}}
\address{Department of Mathematical Analysis, Faculty of Mathematics and Physics, Charles University, Sokolovsk\'a 83, 186 75 Praha 8, Czech Republic}
\email{slavikova@karlin.mff.cuni.cz}

\thanks{B. Park was supported in part by NRF grant 2022R1F1A1063637 and by POSCO Science Fellowship of POSCO TJ Park Foundation. G. Dosidis and L. Slav{\'i}kov{\'a} were supported by the Primus research programme PRIMUS/21/SCI/002 of Charles
University. L. Slav{\'i}kov{\'a} was also partly supported by Charles University Research Centre program No. UNCE/24/SCI/005.}

\title[Bilinear Rough Singular Integrals near the Critical Integrability]{Bilinear Rough Singular Integrals
near the Critical Integrability via Sharp Fourier Multiplier Criteria}
\subjclass[2020]{Primary 42B15, 42B20, 42B25, 47H60}
\keywords{Fourier multiplier operator, Rough singular integral operator, Bilinear operator, Littlewood-Paley theory, Shifted square function}

\begin{abstract}  
We establish  boundedness results for bilinear singular integral operators with rough homogeneous kernels whose restriction to the unit sphere belongs to the Orlicz space $L(\log L)^\alpha$.
This improves the previously best known condition for boundedness of such bilinear operators obtained in the paper of the first and third authors \cite{Do_Sl_submitted}, and provides estimates close to the conjectured endpoint of integrability suggested by the linear theory.
The proof is based on a new sharp boundedness criterion for bilinear Fourier multiplier operators associated with sums of dyadic dilations of a fixed symbol $m_0$, compactly supported away from the origin.
This criterion admits the best possible behavior with respect to a modulation of $m_0$ and is intimately connected with sharp shifted square function estimates.
\end{abstract}

\maketitle

\section{Introduction}\label{S:intro}

With a given bounded measurable function $m$ on $\R^n$ we associate the Fourier multiplier operator 
\begin{equation}\label{E:multiplier}
Tf(x):=\int_{\R^n} m(\xi) \wh{f}(\xi) e^{2\pi i \langle x, \xi\rangle} \,d\xi, \quad x\in \R^n,
\end{equation}
where $f$ is a Schwartz function on $\bbrn$ and $\wh{f}(\xi):=\int_{\bbrn}f(x)e^{-2\pi i\langle x,\xi\rangle}\, dx$ is the Fourier transform of $f$. 
The question of characterizing those symbols $m$ for which the associated Fourier multiplier operator extends to a bounded operator on a certain function space, such as the Lebesgue space $L^p(\R^n)$ or the Hardy space $H^p(\R^n)$, lies at the heart of harmonic analysis. 
In this article, we focus on the class of symbols obtained by summing all dyadic dilations of a given bounded function $m_0$ supported in the unit annulus $\{\xi \in \R^n:~1/2\leq |\xi| \leq 2\}$, namely,
\begin{equation}\label{E:symbol}
m=\sum_{j\in \Z} m_0(2^j\cdot).
\end{equation}
The investigation of this class of symbols is motivated by applications to boundedness of bilinear singular integral operators with rough homogeneous kernels whose restriction to the unit sphere belongs to a function space close to $L^1$. These boundedness results will follow from a new sharp criterion for the boundedness of bilinear Fourier multiplier operators associated with symbols of the form~\eqref{E:symbol}. Since the corresponding Fourier multiplier theorem does not seem to be available in the literature even in the linear setting, we will first briefly discuss the linear case and then move on to the bilinear setting, where the main contributions of our paper lie.


Let $m_0$ and $m$ be as in~\eqref{E:symbol}.
Denoting $K=m_0^{\vee}$, the inverse Fourier transform of $m_0$, we observe that the associated operator $T$, defined in \eqref{E:multiplier}, can be written as
\begin{equation}\label{E:T}
Tf=\sum_{j\in \Z} K_j \ast f.
\end{equation}
 Here and in what follows, $h_j=2^{jn} h(2^{j}\cdot)$ stands for the $L^1$ dyadic dilation of a given function $h$ on   $\bbrn$.

Our point of departure is the following sufficient condition for the boundedness of the operator $T$. It involves the $L^p$-boundedness of $T$ when $1<p<\infty$, as well as the endpoint cases featuring $H^1$ and $BMO$ bounds. To simplify the formulation of the result, we denote
\begin{equation}\label{E:xp}
X^p(\R^n)=\begin{cases}
H^1(\R^n) & \text{if } p=1,\\
L^p(\R^n) &\text{if } 1<p<\infty,\\
BMO(\R^n) & \text{if } p=\infty.
\end{cases}
\end{equation}
As the Schwartz space $\mathscr{S}(\bbrn)$ is not fully contained in $H^1(\bbrn)$, we will study the behavior of the operator $T$ on the subspace $\mathscr{S}_0(\bbrn)$ of the Schwartz class consisting of those $f$ for which
$$\int_{\bbrn}x^{\alpha}f(x)\, dx=0 \q \text{ for all multi-indices }~\alpha.$$
We note that $L^p(\R^n)$ coincides with $H^p(\R^n)$ for any $1<p<\infty$ and that $\mathscr{S}_0(\bbrn)$ is dense in $H^p(\bbrn)$ for all $1\le p<\infty$ (in fact, this holds for every $0<p<\infty$, see~\cite[Chapter III, 5.2(a)]{St1993} or \cite[Theorem 5.1.5]{Tr}).

\begin{proposition}\label{T:Hormander}
Let $1\le p\le \infty$ and $\lambda \geq 0$.
Let $K$ be an integrable function on $\R^n$ whose Fourier transform is supported in the set $\{\xi\in \R^n:~1/2 \leq |\xi| \leq 2\}$ and which satisfies
\begin{equation}\label{E:D}
D_\lambda(K):=\int_{\R^n} |K(y)| \big(\log(e+|y|)\big)^{\lambda}\,dy <\infty.
\end{equation}
Let $T$ be the linear operator associated with $K$, defined by~\eqref{E:T}. 
If $\lambda \ge |1/2-1/p|$,
then there is a constant $C = C(n,p)$, such that
\begin{equation}\label{E:opt-est}
\|Tf\|_{X^p(\R^n)} \le C D_\lambda(K) \|f\|_{X^p(\R^n)}
\end{equation}
for all $f\in\mathscr{S}_0(\bbrn)$.

Conversely, if $\lambda<|1/2-1/p|$, then for each $C>0$ there exists $f\in \mathscr{S}_0(\bbrn)$ and an integrable function $K$ whose Fourier transform is supported in the annulus $\{\xi\in \R^n:~1/2 \leq |\xi| \leq 2\}$ such that
\[
\|Tf\|_{X^p(\R^n)} > C D_\lambda(K)\|f\|_{X^p(\R^n)}.
\]
\end{proposition}

A family of kernels $K$ which is of particular significance in connection with Proposition~\ref{T:Hormander} consists of shifts of a given Schwartz function; this corresponds to the family of modulations of a given Schwartz symbol $m_0$. On the one hand, this family of kernels shows the sharpness of the assumption $\lambda \geq |1/2-1/p|$. At the same time, the family of operators~\eqref{E:T} associated with the above kernels is closely related to the family of shifted square functions, bounds for which play a crucial role in the proof of the sufficiency part of Proposition~\ref{T:Hormander}.

We point out that Proposition~\ref{T:Hormander} is not completely new as various special cases of it have already appeared in the literature. In particular, let us assume that $m$ is a general bounded function, not necessarily having the form~\eqref{E:symbol}, and
instead of \eqref{E:D} assume that
\begin{equation}\label{E:herz-new}
\sup_{k\in \Z} \int_{\R^n} \Big|\big(m(2^k \cdot) \Psi(\cdot)\big)^{\vee}(y)\Big|\big(\log(e+|y|)\big)^\lambda \,dy<\infty,
\end{equation}
where $\Psi$ is a smooth function supported in the set $\{\xi \in \R^n:~1/2 \leq |\xi| \leq 2\}$ and satisfying $\sum_{k\in \Z} \Psi(2^k \xi)=1$ for $\xi \neq 0$.
Then the corresponding multiplier operator $T$, defined in \eqref{E:multiplier}, is bounded from $X^p(\R^n)$ to itself for $1\le p\le \infty$, provided that $\lambda>|1/2-1/p|$. The special case $p=1$ of this result was obtained by Baernstein and Sawyer \cite[Theorem 3b]{BS} and the extension to general $1\le p\le \infty$ was established by Seeger \cite[Theorem 2.1]{Se1990}, who also considered the more general problem of boundedness of the operator $T$ on Triebel-Lizorkin spaces. 
Notably, the strict inequality $\lambda>|1/2-1/p|$ is then required in order to guarantee the $X^p$-boundedness of $T$; for the optimality of this condition, see~\cite[Theorem 4.2]{Park2019_2}. Therefore, the specific form of the symbol~\eqref{E:symbol} is indispensable for the validity of Proposition~\ref{T:Hormander} in the limiting case $\lambda=|1/2-1/p|$.

Related sufficient conditions for the boundedness of Fourier multiplier operators can be found also in~\cite{Car,Park2019_2, Seeger88, Tomita}. Additionally, it is worth noting that the classical H\"ormander multiplier theorem~\cite{Hoe} and its $H^1$-variant due to Calder\'on and Torchinski~\cite{CT} can be recovered as special cases of the results discussed above. More precisely, as pointed out in~\cite[Section 3]{BS}, whenever H\"ormander's condition
\begin{equation*}
\sup_{k\in \Z} \big\|\partial^{\alpha} [m(2^k \cdot) \Psi(\cdot)]\big\|_{L^2(\R^n)} <\infty
\end{equation*}
is satisfied for all multi-indices $\alpha$ with $|\alpha| \leq \lfloor n/2 \rfloor +1$, then
\begin{equation*}
\sup_{k\in \Z} \int_{\R^n} \big|\big(m(2^k \cdot) \Psi(\cdot)\big)\widecheck{}(y)\big|(1+|y|)^\varepsilon \,dy<\infty
\end{equation*}
holds for some $\varepsilon>0$, and, the more so,~\eqref{E:herz-new} is fulfilled.

\hfill

The main result of this article is a sharp bilinear variant of Proposition~\ref{T:Hormander}, stated in Theorem~\ref{T:HormanderB} below, and its application to rough bilinear singular integral operators, given in Theorem~\ref{T:bilinearresult}. The statement of Theorem~\ref{T:HormanderB} involves the bilinear Fourier multiplier operator 
\begin{equation}\label{E:bilinear-fm}
\mathcal{B}(f_1,f_2)(x):=\int_{\R^{2n}} m(\xi_1,\xi_2) \wh{f_1}(\xi_1) \wh{f_2}(\xi_2) e^{2\pi i \langle x,\xi_1+\xi_2\rangle}\,d\xi_1 d\xi_2, \quad x\in \R^n,
\end{equation}
initially defined for $f_1,f_2\in \mathscr{S}_0(\bbrn)$ and associated with a symbol of the form~\eqref{E:symbol}, where $m_0$ is now a function on $\R^{2n}$ supported in the set $\{(\xi_1,\xi_2)\in \R^{2n}:~1/2 \leq |(\xi_1,\xi_2)| \leq 2\}$. Denoting again $K=m_0^{\vee}$, we observe that $\mathcal {B}$ can be rewritten as
\begin{equation}\label{E:b}
\mathcal{B}(f_1,f_2)(x)=\sum_{l\in \Z} K_{l} * \big(f_1 \otimes f_2\big) (x,x).
\end{equation}
We recall that $K_l = 2^{2ln} K(2^{l} \cdot)$. We will establish a sufficient condition for the $L^{p_1}(\R^n) \times L^{p_2}(\R^n) \rightarrow L^p(\R^n)$ boundedness of $\mathcal{B}$ whenever $1<p_1, p_2 \leq \infty$, $1 \leq p<\infty$ and the exponents satisfy the H\"older scaling $\frac{1}{p_1}+\frac{1}{p_2}=\frac{1}{p}$. We will also prove endpoint bounds involving $H^1(\R^n)$ as one of the input spaces, or $BMO(\R^n)$ as the target space. Using the fact that the Lebesgue space $L^p(\R^n)$ coincides with the Hardy space $H^p(\R^n)$ for $1<p\leq \infty$, all input spaces will be Hardy spaces in the next theorem. We also introduce the scale of target spaces $Y^p(\R^n)$ as
$$
Y^p(\R^n)=\begin{cases}  L^p(\R^n) & \text{if } 1\le p<\infty,\\
BMO(\R^n) & \text{if } p=\infty.  \end{cases}
$$
Given $p\in [1,\infty]$, we employ the usual convention that  $p'\in [1,\infty]$ stands for the conjugate exponent to $p$ satisfying $\frac{1}{p}+\frac{1}{p'}=1$.

\begin{theorem}\label{T:HormanderB}
Let $1\le p,p_1,p_2\le \infty$ and $\frac1p= \frac1{p_1} + \frac{1}{p_2}$.
Suppose that $K$ is an integrable function on $\R^{2n}$ whose Fourier transform is supported in the set $\{(\xi_1,\xi_2)\in \R^{2n}:~1/2 \leq |(\xi_1,\xi_2)| \leq 2\}$ and which satisfies
\begin{equation*}
D_\lambda(K):=\int_{\R^{2n}} \big|K(\yyy)\big| \big(\log(e+|\yyy|\big))^{\lambda}\,d\yyy <\infty
\end{equation*}
where  $\yyy=(y_1,y_2)\in \R^{2n}$.
Let $\mathcal{B}$ be the bilinear operator associated with $K$, defined by~\eqref{E:b}.
If
\begin{equation}\label{E:LambdaB} \lambda\geq\max\Big\{\frac1{p_1}, \frac{1}{p_2}, \frac{1}{p'}\Big\},
\end{equation}
then there is a constant $C= C(n,p_1,p_2)$ such that
\begin{equation}\label{E:BoundB} 
\|\mathcal{B}(f_1,f_2)\|_{Y^p(\R^n)} \leq C D_\lambda(K) \|f_1\|_{H^{p_1}(\R^n)}\|f_2\|_{H^{p_2}(\R^n)}
\end{equation}
for $f_1,f_2\in\mathscr{S}_0(\bbrn)$.

Conversely, if $\lambda<\max\Big\{\frac1{p_1}, \frac{1}{p_2}, \frac{1}{p'}\Big\}$, then for each $C>0$ there exists a pair of functions $f_1, f_2\in\mathscr{S}_0(\bbrn)$ and an integrable function $K$ whose Fourier transform is supported in the annulus  $\{(\xi_1,\xi_2)\in \R^{2n}:~1/2 \leq |(\xi_1,\xi_2)| \leq 2\}$ such that 
\begin{equation*} 
\big\|\mathcal{B}(f_1,f_2)\big\|_{Y^p(\R^n)} > C D_\lambda(K) \|f_1\|_{H^{p_1}(\R^n)}\|f_2\|_{H^{p_2}(\R^n)}.
\end{equation*}
\end{theorem} 

\begin{remark}\label{remark1}
Theorem \ref{T:HormanderB} is also valid with $f_1,f_2\in \mathscr{S}_0(\bbrn)$ replaced with $f_1,f_2\in\mathscr{S}(\bbrn)$ when $1<p_1,p_2\le \infty$. This setting will be needed in the proof of Theorem \ref{T:bilinearresult} below.
\end{remark}

\begin{figure}[h]
\centering
\scalebox{1.5}{\includegraphics[page=1]{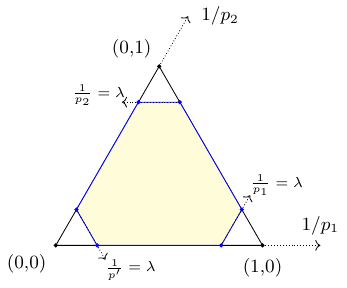}}
\caption{Range of $p_1,p_2$ for $H^{p_1}\times H^{p_2}\to Y^p$ boundedness of $\mathcal{B}$ for a given value of $\lambda$.}
\label{fig:B}
\end{figure}

The proof of Theorem~\ref{T:HormanderB} is based on sharp estimates for the family of shifted square functions combined with delicate symmetry and duality considerations. 

We note that the quantity on the right-hand side of~\eqref{E:LambdaB} is bounded from below by $1/3$, and thus the mere integrability of $K$ is not sufficient to guarantee the $L^{p_1}(\R^n) \times L^{p_2}(\R^n) \rightarrow L^p(\R^n)$ boundedness of the associated bilinear operator for any choice of exponents $p_1$, $p_2$, $p$. This is in stark contrast with the $L^2$ estimate for linear operators in Proposition \ref{T:Hormander}.

\hfill

As an application of Theorem~\ref{T:HormanderB}, we obtain improved estimates for bilinear singular integral operators with rough homogeneous kernels. Each operator of this type is associated with an  integrable function $\Omega$
on $\mathbb{S}^{2n-1}$ satisfying the vanishing moment condition
\begin{equation}\label{E:vanishingcon}
\int_{\bbs^{2n-1}}\Omega(\theta)\;\ d\nu(\theta)=0,
\end{equation}
where $d\nu$ stands for the normalized surface measure on the unit sphere $\bbs^{2n-1}$,
and it acts on a pair of Schwartz functions $f_1$ and $f_2$ as
\begin{equation}\label{E:rough}
T_{\Omega}\big(f_1,f_2\big)(x)=\mathrm{p.v.} \int_{\R^{2n}}{\frac{\Omega(\theta)}{|\yyy|^{2n}} f_{1}(x-y_{1}) f_{2}(x-y_{2})}~d\,\yyy, \quad x\in \R^n,
\end{equation} 
where  $\yyy=(y_1,y_2)\in \R^{2n}$ and $\theta=\yyy/|\yyy\,|\in \mathbb{S}^{2n-1}$.

Boundedness properties of the operator $T_\Omega$ when $\Omega \in L^q(\S^{2n-1})$ for some $q>1$ were studied in~\cite{Gr_He_Ho2018, GHS20, He_Park2022}; see also~\cite{Do_Sl_submitted, Gr_He_Ho_Park_submitted} for results that apply in the more general multilinear setting. It follows from the results of these papers that $T_\Omega$ is $L^{p_1}(\R^n) \times L^{p_2}(\R^n) \rightarrow L^p(\R^n)$ bounded for all $\Omega \in L^q(\S^{2n-1})$  satisfying~\eqref{E:vanishingcon} if and only if 
\begin{equation}\label{E:pq}
1/p+1/q<2.
\end{equation}
Here, $1<p_1, p_2<\infty$ and $1/p=1/p_1+1/p_2$; in particular, $p$ may be smaller than $1$.
In this article, we focus on the limiting case when $\Omega$ is not necessarily $L^q(\S^{2n-1})$ integrable for any $q>1$ but it satisfies
\begin{equation}\label{E:omega}
\int_{\S^{2n-1}} |\Omega(\theta)| \big(\log(e+ |\Omega(\theta)|)\big)^\alpha \, d\nu(\theta)<\infty
\end{equation}
for some $\alpha \geq 0$. If condition~\eqref{E:omega} is fulfilled, we say that $\Omega \in L(\log L)^\alpha(\S^{2n-1})$. The norm on the space $L(\log L)^\alpha(\S^{2n-1})$ is then defined as
\[
\|\Omega\|_{L(\log L)^\alpha(\S^{2n-1})}
=\inf\left\{\lambda>0:~\int_{\S^{2n-1}} \frac{|\Omega(\theta)|}{\lambda} \left[\log\left(e+\frac{|\Omega(\theta)|}{\lambda}\right)\right]^\alpha\,d\nu(\theta) \leq 1\right\}.
\]

As any function $\Omega \in \bigcup_{q>1} L^q(\S^{2n-1})$ fulfills ~\eqref{E:omega} for every $\alpha \geq 0$, the condition~\eqref{E:pq} implies that the operator $T_\Omega$ can only be $L^{p_1}(\R^n) \times L^{p_2}(\R^n) \rightarrow L^p(\R^n)$ bounded for all $\Omega \in L(\log L)^\alpha(\S^{2n-1})$ if $p \geq 1$. We show that the opposite implication is also true, as long as $\alpha \geq 2$. In addition, for values of $\alpha$ from the interval $[4/3,2)$ we obtain the $L^{p_1}(\R^n) \times L^{p_2}(\R^n) \rightarrow L^p(\R^n)$ boundedness of $T_\Omega$ in a restricted range of exponents. Both these claims follow as a consequence of the following theorem. 

\begin{theorem}\label{T:bilinearresult}
Let $1\le p<\infty$ and $1<p_1,p_2\le \infty$ with $\frac{1}{p_1}+\frac{1}{p_2}=\frac{1}{p}$. Suppose that
\begin{equation}\label{E:A}
  A \geq \max\Big\{\frac1{p_1}, \frac{1}{p_2}, \frac{1}{p'}\Big\}+1
\end{equation}
  and 
  $$\Omega\in L(\log L)^A(\S^{2n-1}) \q \text{with}\q\int_{\bbs^{2n-1}} \Omega(\theta)\; d\nu(\theta) = 0.$$
Then there exists a constant $C=C(n,p_1,p_2)$ such that
\[ \big\|T_{\Omega}(f_1,f_2)\big\|_{L^p(\R^n)}\le C \|\Omega\|_{L(\log L)^A(\S^{2n-1})}  \|f_{1}\|_{L^{p_{1}}(\R^n)}  \|f_{2}\|_{L^{p_{2}}(\R^n)} \]
for all Schwartz functions $f_1$ and $f_2$ on $\bbrn$.
\end{theorem}

The proof of Theorem~\ref{T:bilinearresult} employs the approach of Duoandikoetxea and Rubio de Francia~\cite{Du_Ru1986}, decomposing the rough kernel into an infinite sum of smooth kernels. The operators associated with the smooth kernels are then estimated by making use of Theorem~\ref{T:HormanderB} in combination with earlier estimates from~\cite{Gr_He_Ho2018, He_Park2022}.

\medskip

It remains still open whether the exponent $A$ in \eqref{E:A} is sharp.
A comparison with the linear theory
\cite{Ca_Zy1956, Ch_Ru1988, Ho1988, Se1996}
indicates that this may not be the case.
Indeed, it is well known that a linear singular integral operator associated with a
homogeneous kernel whose restriction $\Omega$ to the unit sphere
$\mathbb S^{n-1}$ belongs to the space $L\log L(\mathbb S^{n-1})$ is $L^p$-bounded
for any $1<p<\infty$.
However, such a conclusion (corresponding to $A=1$) cannot be recovered via the
methods of the present article.
Nevertheless, classical methods for proving $L^p$-boundedness of linear rough singular
integrals with $\Omega \in L\log L(\mathbb S^{n-1})$, such as the method of rotations,
do not extend easily to the bilinear setting.
This limitation motivates the alternative approach developed in the present paper.
To the best of our knowledge, our work is the first to investigate the applicability of
the techniques of Duoandikoetxea and Rubio de Francia~\cite{Du_Ru1986} within the
$L(\log L)^\alpha$ framework, providing a flexible methodology that extends naturally
to bilinear and even multilinear operators.
While the sharp results in Proposition~\ref{T:Hormander} and
Theorem~\ref{T:HormanderB} reveal that this approach may not suffice to reach the full
endpoint condition $\Omega \in L\log L$, it does yield a nearly optimal result,
namely $\Omega \in L(\log L)^2$, in the limiting integrability case $q=1$.
Moreover, an important feature of the present approach is that this phenomenon is not
confined to the bilinear setting.
Indeed, the same mechanism persists in the multilinear case, as demonstrated in the
recent preprint~\cite{Haar}, highlighting both the robustness of the method and its
potential relevance for the development of a broader endpoint theory for multilinear
operators with rough kernels.

\medskip

{\bf Organization of the paper.}
Section~\ref{S:preliminary} contains some preliminary materials involving the Hardy spaces $H^p(\R^n)$, the space $BMO(\R^n)$, and certain maximal inequalities. In Section~\ref{S:shifted}, we obtain sharp shifted square function estimates, which will provide a key technical ingredient for the proofs of our main results. 
The proof of Proposition~\ref{T:Hormander} is given in Section~\ref{S:linear}, and the proofs of Theorems \ref{T:HormanderB} and \ref{T:bilinearresult} are the content of Sections \ref{S:bilinear}  and \ref{proof1}, respectively.
\\

{\bf Notation.}
We denote by $\bbn$ and $\bbz$ the sets of   natural numbers and   integers, respectively, and let $\bbn_0:=\bbn\cup \{0\}$. We use the symbol $A\lesssim B$ to indicate that $A\leq CB$ for some constant $C>0$ independent of the variable quantities $A$ and $B$, and $A\sim B$ if $A\lesssim B$ and $B\lesssim A$ hold simultaneously.
For each cube $Q$ in $\bbrn$, we denote by $\ell(Q)$ the side length of $Q$ and by $Q^*$ the cube with the same center as $Q$ and $10\sqrt{n}$ times the side length.
The set of all dyadic cubes in $\bbrn$ is denoted by   $\mathcal{D}$, and for each $k\in\mathbb{Z}$ we 
designate    $\mathcal{D}_{k}$ to be  the subset of  $\mathcal{D}$ consisting of dyadic cubes with side length $2^{-k}$. 

\section{Preliminaries}\label{S:preliminary}

Throughout the paper, $\phi$ and $\psi$ will stand for Schwartz functions on $\bbrn$ such that
$$\wh{\phi}(0)=1,\q \supp{(\wh{\phi})}\subset \{\xi\in\bbrn: |\xi|\lesssim 1\},$$
\begin{equation}\label{psiprop}
\supp{(\wh{\psi})}\subset \{\xi\in\bbrn: |\xi|\sim 1\}, \q \text{ and }\q \sum_{k\in\bbz}\wh{\psi}(2^k\xi)=C,~ \xi\not= 0
\end{equation}
for some $C>0$.

\subsection{Hardy spaces and $BMO$}
 For $0<p\leq \infty$, the Hardy space $H^p(\bbrn)$ consists of all tempered distributions $f$ on $\bbrn$ satisfying
\begin{equation*}
\Vert f\Vert_{H^p(\bbrn)}:=\Big\Vert \sup_{k\in\bbz}{|\phi_k\ast f|} \Big\Vert_{L^p(\bbrn)}<\infty.
\end{equation*}
The Hardy space coincides with the Lebesgue space $L^p(\bbrn)$ for $1<p\le \infty$, while $H^p(\R^n)\subsetneq L^p(\R^n)$ for $0<p\le 1$.
According to \cite{Fr_Ja1990}, different choices of the Schwartz function $\phi$ provide equivalent Hardy space norms.
Moreover, the space $H^p(\R^n)$ can be characterized in terms of 
 Littlewood-Paley theory. For $0<p<\infty$  we have
\begin{equation}\label{hardycha}
\Vert f\Vert_{H^p(\bbrn)}\sim \big\Vert \{ \psi_k\ast f\}_{k\in\bbz} \big\Vert_{L^{p}(\ell^2)}.
\end{equation}

For $p=\infty$, the equivalence \eqref{hardycha} does not hold, but we have the following characterization of the $BMO$ space, which originated from a Carleson measure characterization:
\begin{equation}\label{bmocha}
\Vert f\Vert_{BMO(\bbrn)} \sim \sup_{P\in\mathcal{D}}\bigg( \frac{1}{|P|}\int_P\sum_{k=-\log_2{\ell(P)}}^{\infty}\big| \psi_k\ast f(x)\big|^2\;dx\bigg)^{1/2}.
\end{equation}
Here, the supremum is taken over all dyadic cubes $P$.
The properties \eqref{hardycha} and \eqref{bmocha} are also independent of the choice of functions $\psi_k$ (see, e.g. \cite[Remark 2.6, Corollary 5.3]{Fr_Ja1990}).

\subsection{Maximal inequalities}\label{maximalines}

Let $1<r<\infty$. We denote by $\MM_r$ the modified Hardy-Littlewood maximal operator defined by
$$\MM_r f(x):=\sup_{Q:x\in Q}\left(\frac{1}{|Q|}\int_Q{|f(y)|^r}dy\right)^{\frac{1}{r}},$$
where the supremum is taken over all cubes $Q$ containing $x$. Given $k\in \Z$ and $\sigma>0$, we also introduce Peetre's maximal function
\begin{equation*}
\mathfrak{M}_{\sigma,2^k}f(x):=\sup_{y\in\bbrn}{\frac{|f(x-y)|}{(1+2^k|y|)^{\sigma}}}.
\end{equation*}
 For $r>0$, let $\mathcal{E}(r)$ denote the space of all distributions whose Fourier transforms are supported in $\big\{\xi\in\bbrn:|\xi|\leq 2r\big\}$.
 It was shown in \cite{Pe1975} that
 \begin{eqnarray}\label{maximalbound}
\mathfrak{M}_{n/r,2^k}f(x)\lesssim_{r,A} \mathcal{M}_rf(x), 
\end{eqnarray}  provided that $f\in\mathcal{E}(A2^k)$ for $A>0$.
A combination of~\eqref{maximalbound} and of the Fefferman-Stein vector-valued maximal inequality~\cite{Fe_St1971} yields that 
 for $0<p<\infty$ and $0<q\leq \infty$, we have 
\begin{equation}\label{peetremax}
\big\Vert  \big\{\mathfrak{M}_{\sigma,2^k}f_k \big\}_{k\in\bbz} \big\Vert_{L^p(\ell^q)} \lesssim_{A,p,q}  \big\Vert \big\{f_k\big\}    \big\Vert_{L^p(\ell^q)} \quad \text{for}~\sigma>n/\min{(p,q)}
\end{equation}  if $f_k\in\mathcal{E}(A2^k)$.
Clearly, the above inequality also holds for $p=q=\infty$.

A substitute for the inequality~\eqref{peetremax} in the case when $p=\infty$ and $0<q<\infty$ has the form
\begin{equation}\label{bmomaximaline}
\sup_{P\in\mathcal{D}}{\bigg( \frac{1}{|P|}\int_P{\sum_{k=-\log_2{\ell(P)}}^{\infty}{\big(\mathfrak{M}_{\sigma,2^k} f_k(x) \big)^q}}\;dx\bigg)^{1/q}}\lesssim_{A,q} \sup_{P\in\mathcal{D}}{\bigg( \frac{1}{|P|}\int_P{\sum_{k=-\log_2{\ell(P)}}^{\infty}{\big|f_k(x) \big|^q}}\;dx\bigg)^{1/q}}
\end{equation} for $\sigma>n/q$ if $f_k\in\mathcal{E}(A2^k)$.
We note that the inequality \eqref{bmomaximaline} with $\mathfrak{M}_{\sigma,2^k}$ replaced by $\MM_r$ is false for any $r>0$.
 See \cite{Park2019} for more details.

\subsection{A vector-valued characterization of $H^p$ and $BMO$}
 Let $\sigma>0$ and $k\in \Z$. It is straightforward to verify that the function $\mathfrak{M}_{\sigma,2^k}$ is essentially constant on each $Q\in \mathcal{D}_k$, namely
\begin{equation}\label{infmax2}
\sup_{y\in Q}{\mathfrak{M}_{\sigma,2^k}f(y)} \lesssim \inf_{y\in Q}{\mathfrak{M}_{\sigma,2^k}f(y)}.
\end{equation}
Moreover, we have the following result which, when combined with~\eqref{hardycha} and~\eqref{bmocha}, yields a vector-valued characterization of $H^p$ and $BMO$ spaces.
\begin{customlemma}{A}[\cite{Park_IUMJ}]\label{equivalencelemma}
Let $0<p\le \infty$, $0<\gamma<1$, and $\sigma>n/\min{(p,2)}$.
Suppose that $\fff:=\{f_k\}_{k\in\bbz}$ is a sequence of $f_k\in \mathcal{E}(A2^k)$ for some $A>0$.
Then for each dyadic cube $Q\in\mathcal{D}$, there exists a proper measurable subset $S_Q$ of $Q$, depending on $\gamma,\fff$, such that $|S_Q|>\gamma |Q|$ and 
\begin{equation}\label{finitepsqureest}
 \bigg\Vert \Big\{ \sum_{Q\in\mathcal{D}_k}\Big( \inf_{y\in Q}\mathfrak{M}_{\sigma,2^k}{f_k}(y)\Big)\chi_{S_Q}\Big\}_{k\in\bbz}\bigg\Vert_{L^p(\ell^2)}\sim \big\Vert \fff \big\Vert_{L^p(\ell^2)} \quad \text{ if }~0<p<\infty
 \end{equation}
and
\begin{align*}
&\bigg\Vert \Big\{ \sum_{Q\in\mathcal{D}_k}\Big( \inf_{y\in Q}\mathfrak{M}_{\sigma,2^k}{f_k}(y)\Big)\chi_{S_Q}\Big\}_{k\in\bbz}\bigg\Vert_{L^{\infty}(\ell^2)}\\
&\sim \sup_{P\in\mathcal{D}}\Big(\frac{1}{|P|}\int_P \sum_{k=-\log_2{\ell(P)}}^{\infty}\big| f_k(x)\big|^2\; dx \Big)^{1/2} \quad \text{ if }~p=\infty.
\end{align*}
\end{customlemma}
Note that the constants in both equivalences are independent of $\fff$, while the dependence on $\fff$ is necessary in choosing a subset $S_Q$.


\section{Shifted square function estimates}\label{S:shifted}
Let $\psi$ be a Schwartz function satisfying \eqref{psiprop}.
For $k\in\bbz$ and $y\in\bbrn$, we define
$$(\psi_k)^y:=\psi_k(\cdot-2^{-k}y)=2^{kn}\psi(2^k\cdot-y).$$
The main goal of this section is to prove the following shifted square function estimates, which will become one of the key ingredients in the proof of Proposition~\ref{T:Hormander} and Theorem~\ref{T:HormanderB}.

\begin{lemma}\label{lpbmoest}
Let $y\in\bbrn$.
Then for $1 \leq p< \infty$
\begin{equation}\label{squarefiniteest}
\Big\Vert \Big( \sum_{k\in\bbz}\big|(\psi_k)^y\ast f\big|^2 \Big)^{1/2}\Big\Vert_{L^p(\bbrn)}\lesssim \big(\log{(e+|y|)}\big)^{|1/p-1/2|}\Vert f \Vert_{H^p(\bbrn)},
\end{equation}
and for $p=\infty$
\begin{equation}\label{suppdbmoest}
\sup_{P\in\mathcal{D}}\bigg(\frac{1}{|P|}\int_P\sum_{k=-\log_2{\ell(P)}}^{\infty}\big| (\psi_k)^y\ast f(x)\big|^2\; dx \bigg)^{1/2}\lesssim \big( \log{(e+|y|)}\big)^{1/2}\Vert f\Vert_{BMO(\bbrn)}.
\end{equation}
\end{lemma}

Variants of the estimate~\eqref{squarefiniteest} with the exponent $|1/p-1/2|$ replaced by $1$ are available in the literature, see~\cite[Theorem 5.1]{Mu2014} in the one-dimensional setting and~\cite[Proposition 7.5.1]{MFA} for the higher-dimensional variant. 
The logarithmic bounds in~\eqref{squarefiniteest} and~\eqref{suppdbmoest} are sharp, in the sense that the exponent $|1/p-1/2|$ cannot be replaced by any lower number. This follows by an argument analogous to the one given in Section~\ref{SS:OptimalityLinear} below.


A combination of Lemmas~\ref{equivalencelemma} and~\ref{lpbmoest} then yields the following corollary, which provides a unified treatment of the shifted square function estimates for $H^p(\R^n)$ and $BMO(\R^n)$. We recall that $X^p(\R^n)$ stands for the scale of function spaces defined in~\eqref{E:xp}.

\begin{corollary}\label{C:ShiftedSquare}
Let $1\le p\le \infty$, $0<\gamma<1$ and $y\in\bbrn$.
Suppose that $f\in X^p(\R^n)$ and $\sigma>n/\min{(p,2)}$.
Then for each dyadic cube $Q\in\mathcal{D}$, there exists a proper measurable subset $S_Q$ of $Q$, depending on $\gamma,y$ and $f$, such that $|S_Q|>\gamma |Q|$ and
\begin{equation*}
\bigg\Vert \Big\{ \sum_{Q\in\mathcal{D}_k}\Big( \inf_{z\in Q}\mathfrak{M}_{\sigma,2^k}{\big((\psi_k)^y\ast f\big)}(z)\Big)\chi_{S_Q}\Big\}_{k\in\bbz}\bigg\Vert_{L^p(\ell^2)}\lesssim \big(\log{(e+|y|)} \big)^{|1/p-1/2|}\Vert f \Vert_{X^p(\bbrn)}.
\end{equation*}
\end{corollary}


We next focus on proving Lemma~\ref{lpbmoest}. We start with
the following auxiliary result, which is the key estimate for establishing inequality \eqref{suppdbmoest}.

\begin{lemma}\label{linfkeyest}
Let $y\in\bbrn$ and $P\in\mathcal{D}$. Suppose that each $f_k$ belongs to $\mathcal{E}(A2^k)$ for some $A>0$.
Then we have
\begin{align}\label{E:log}
&\bigg(\frac{1}{|P|}\int_P   \sum_{k=-\log_2{\ell(P)}}^{\infty}\big| (\psi_k)^y \ast f_k(x)\big|^2     \, dx \bigg)^{1/2}\nonumber\\
&\lesssim (\log{(e+|y|)})^{{1}/{2}}\sup_{R\in\mathcal{D}:\ell(R)\le \ell(P)}\bigg(\frac{1}{|R|}\int_R   \sum_{k=-\log_2{\ell(R)}}^{\infty}\big|  f_k(x)\big|^2     \, dx \bigg)^{1/2}.
\end{align}
\end{lemma}
\begin{proof}
We bound the left-hand side of~\eqref{E:log} by
\begin{align*}
&\bigg(\frac{1}{|P|}\int_P   \sum_{k=-\log_2{\ell(P)}}^{\infty}\big| (\psi_k)^y \ast \big( f_k \cdot \chi_{P^*}  \big)(x)\big|^2     \, dx \bigg)^{1/2}\\
&\quad + \bigg(\frac{1}{|P|}\int_P   \sum_{k=-\log_2{\ell(P)}}^{\infty}\big| (\psi_k)^y \ast \big( f_k\cdot \chi_{(P^{*})^c}  \big)(x)\big|^2     \, dx \bigg)^{1/2}\\
&=:\mathrm{I}_P+\mathrm{II}_P.
\end{align*}
We recall that $P^*$ denotes the concentric dilate of $P$ with $\ell(P^*)=10\sqrt{n}\ell(P)$.

The term $\mathrm{I}_P$ is bounded by
\begin{equation}\label{I'Pest}
\bigg(\frac{1}{|P|}\sum_{k=-\log_2{\ell(P)}}^{\infty}\big\Vert (\psi_k)^y\ast \big(f_k\cdot \chi_{P^*} \big)\big\Vert_{L^2(\bbrn)}^2 \bigg)^{1/2}
\end{equation}
and Young's inequality for convolutions gives
\begin{align*}
\big\Vert (\psi_k)^y\ast \big(f_k\cdot \chi_{P^*} \big)\big\Vert_{L^2(\bbrn)}^2\le \big\Vert (\psi_k)^y\big\Vert_{L^1(\bbrn)}^2\big\Vert f_k\big\Vert_{L^2(P^*)}^2\lesssim \int_{P^*}\big| f_k(x)\big|^2 \; dx,
\end{align*}
which implies that \eqref{I'Pest} is controlled by a constant times
\begin{align*}
&\bigg(\frac{1}{|P|}\int_{P^*}\sum_{k=-\log_2{\ell(P)}}^{\infty}\big| f_k(x)\big|^2 \; dx\bigg)^{1/2}\\
&\lesssim \sup_{R\in\mathcal{D}:~\ell(R)=\ell(P)}\bigg( \frac{1}{|R|}\int_{R}\sum_{k=-\log_2{\ell(R)}}^{\infty}\big| f_k(x)\big|^2 \; dx\bigg)^{1/2}.
\end{align*}
This proves the desired estimate for $I_P$.

Now we consider the term $\mathrm{II}_P$, for which the constant $\big( \log{(e+|y|)}\big)^{1/2}$ is required in the upper bound.
We write
\begin{align*}
\mathrm{II}_P&\le \bigg(\frac{1}{|P|}\int_P  \sum_{k:~1\le 2^k\ell(P)<1+|y|}  \big| (\psi_k)^y \ast \big( f_k \cdot \chi_{(P^{*})^c}  \big)(x)\big|^2    \; dx \bigg)^{1/2}\\
& \qq+ \bigg(\frac{1}{|P|}\int_P  \sum_{k:~2^k\ell(P)\ge 1+|y|}  \big| (\psi_k)^y \ast \big( f_k \cdot \chi_{(P^{*})^c}  \big)(x)\big|^2    \; dx \bigg)^{1/2}\\
&=:\mathrm{II}_P^1+\mathrm{II}_P^2.
\end{align*}
Given $x\in P$ and $k\in \Z$, we have
\[
\big| (\psi_k)^y \ast \big( f_k \cdot    \chi_{(P^{*})^c}  \big)(x)\big|
\lesssim \|(\psi_k)^y\|_{L^1(\R^n)}\|f_k\|_{L^\infty(\R^n)}
\lesssim \|f_k\|_{L^\infty(\R^n)}.
\]
It was verified in \cite[Section 1]{Park2019} that
\begin{equation}\label{fkinfest}
\Vert  f_k\Vert_{L^{\infty}(\bbrn)}\lesssim_A \sup_{R\in\mathcal{D}:2^k\ell(R)\le 1}\bigg(\frac{1}{|R|}\int_R\sum_{j=-\log_2{\ell(R)}}^{\infty}\big|  f_j(x)\big|^2\; dx \bigg)^{1/2},
\end{equation}
provided that each $f_j\in\mathcal{E}(A2^j)$ for some $A>0$.

This yields
\begin{align*}
\mathrm{II}_P^1
\lesssim \big(\log{(e+|y|)}\big)^{1/2} \sup_{R\in\mathcal{D}:\ell(R)\le \ell(P)}\bigg(\frac{1}{|R|}\int_R\sum_{j=-\log_2{\ell(R)}}^{\infty}\big|  f_j(x)\big|^2\; dx \bigg)^{1/2}.
\end{align*}
To estimate the term $\mathrm{II}_P^2$, we observe that if $2^k\ell(P)\ge 1+|y|$, then $|2^{-k}y|\le\ell(P)$, and thus for $x\in P$ and $z\in (P^*)^c$, we have
$$|x-z-2^{-k}y|\ge |x-z|-|2^{-k}y|>|x-z|-\ell(P)\gtrsim |x-z|\gtrsim \ell(P).$$
This implies that
\begin{align*}
& \big| (\psi_k)^y \ast \big( f_k \cdot    \chi_{(P^{*})^c}  \big)(x)\big|\\
&\lesssim \int_{(P^*)^c}  \frac{2^{kn}}{\big( 1+2^k|x-z-2^{-k}y|\big)^{n+2}} \big| f_k(z)\big|        \; dz\\
&\lesssim   \big( 2^k\ell(P)\big)^{-1} \Vert f_k \Vert_{L^{\infty}(\bbrn)}  \int_{(P^*)^c}\frac{2^{kn}}{\big(1+2^k|x-z|\big)^{n+1}}\; dz\\
&\lesssim \big( 2^k\ell(P)\big)^{-1} \Vert f_k\Vert_{L^{\infty}(\bbrn)}.
\end{align*}
Summing a geometric series and using \eqref{fkinfest} again, we thus deduce
\begin{align*}
\mathrm{II}_P^2
\lesssim \sup_{R\in\mathcal{D}:\ell(R)\le \ell(P)}\bigg(\frac{1}{|R|}\int_R\sum_{j=-\log_2{\ell(R)}}^{\infty}\big|  f_j(x)\big|^2\; dx \bigg)^{1/2},
\end{align*}
which completes the proof of \eqref{E:log}.
\end{proof}

Now let us prove Lemma \ref{lpbmoest}.

\begin{proof}[Proof of Lemma \ref{lpbmoest}]
We first verify \eqref{suppdbmoest}.
Thanks to \eqref{psiprop}, there exists a positive integer $k_0$ such that
\begin{equation}\label{E:psi-tilde}
\wh{\wt{\psi}}(\xi):=\frac{1}{C}\sum_{j=-k_0}^{k_0}\wh{\psi_j}(\xi)= 1 \q \text{ on }~ \supp{(\wh{\psi})}.
\end{equation}
Note that condition~\eqref{psiprop} is then satisfied with $\psi$ replaced by $\wt{\psi}$ and with $C=2k_0+1$. In addition, \eqref{E:psi-tilde} implies
\begin{equation}\label{simest}
 (\psi_k)^y\ast f=  (\psi_k)^y\ast \big(\wt{\psi}_k\ast f \big)
\end{equation}
for $k\in \mathbb Z$.
Using \eqref{simest} and Lemma \ref{linfkeyest} with $f_k=\wt{\psi}_k\ast f$, the left-hand side of \eqref{suppdbmoest} is controlled by a constant times
\begin{equation*}
\big( \log{(e+|y|)}\big)^{1/2}\sup_{P\in\mathcal{D}}\bigg(\frac{1}{|P|}\int_P\sum_{k=-\log_2{\ell(P)}}^{\infty}\big| \wt{\psi}_k\ast f(x)\big|^2\; dx \bigg)^{1/2},
\end{equation*}
which is comparable to the right-hand side of \eqref{suppdbmoest} due to \eqref{bmocha}. This proves \eqref{suppdbmoest}.

Now we consider the case when $1\le p<\infty$.
We choose an integer $\rho_0$ such that
$$\supp(\wh{\psi})\subset \big\{\xi\in\bbrn : 2^{-\rho_0}\le |\xi|\le 2^{\rho_0}\big\}.$$
We split the set $\mathbb Z$ into $10(\rho_0+k_0)$ subsets modulo $10(\rho_0+k_0)$. We fix one such subset $A$ and define
$$T^{y}f(x):=\sum_{k\in A}(\psi_{k})^y\ast f(x).$$
We first observe that 
\begin{equation}\label{equiclaim}
\big\Vert T^{y}f\big\Vert_{H^p(\bbrn)}\sim_{\rho_0,k_0}\bigg\Vert \Big( \sum_{k\in A}\big| \big(\psi_{k})^y\ast f\big|^2\Big)^{1/2}\bigg\Vert_{L^p(\bbrn)}.
\end{equation}

Indeed, since
 $\psi_l\ast (\psi_{k})^y\ast f(x)$ vanishes unless $k-2\rho_0\le l \le k+2\rho_0$, 
 the left-hand side of \eqref{equiclaim} is comparable,  via \eqref{hardycha}, to
 \begin{equation}\label{sumksuml}
 \big\Vert \big\{ \psi_l\ast T^{y}f\big\}_{l\in\bbz}\big\Vert_{L^p(\ell^2)} = \bigg\Vert \bigg(\sum_{k\in A}\sum_{l=k-2\rho_0}^{k+2\rho_0}\big| \psi_l\ast (\psi_{k})^y\ast f\big|^2 \bigg)^{1/2}\bigg\Vert_{L^p(\bbrn)}.
 \end{equation}
For each $k-2\rho_0\le l \le k+2\rho_0$, we have
\begin{equation}\label{doconest}
\big| \psi_l\ast (\psi_{k})^y\ast f(x) \big|\lesssim_{\sigma,\rho_0,k_0} \mathfrak{M}_{\sigma,2^{k}}\big((\psi_{k})^y\ast f\big)(x) \q \text{ for }~\sigma>n,
\end{equation}
and thus the right-hand side of \eqref{sumksuml} is bounded by a constant times
$$\big\Vert \big\{ \mathfrak{M}_{\sigma,2^{k}}\big( (\psi_{k})^y\ast f\big) \big\}_{k\in\bbz}\big\Vert_{L^p(\ell^2)}.$$
Finally, the maximal inequality \eqref{peetremax} proves one direction of \eqref{equiclaim}.

To verify the opposite direction, we observe that for any $k\in A$,
\[
(\psi_k)^y \ast f= \wt{\psi}_k \ast (\psi_k)^y \ast f
=\wt{\psi}_k \ast T^y f,
\]
where the last equality holds since the elements of the set $A$ are well separated. The right-hand side of~\eqref{equiclaim} is thus equal to
$$\bigg\Vert \bigg(\sum_{k\in A}   \big| \wt{\psi}_k\ast T^{y}f\big|^2        \bigg)^{1/2}\bigg\Vert_{L^p(\bbrn)}\lesssim \big\Vert T^{y}f\big\Vert_{H^p(\bbrn)},$$
as desired.

In view of~\eqref{equiclaim}, inequality~\eqref{squarefiniteest} will follow if we prove
\begin{equation}\label{sufficeest}
\big\Vert T^{y}f\big\Vert_{H^p(\bbrn)}\lesssim (\log(e+|y|))^{|1/p-1/2|}\Vert f\Vert_{H^p(\bbrn)} \q \text{ for }~1\le p<\infty.
\end{equation}
Indeed, once~\eqref{sufficeest} is proved for all $10(\rho_0+k_0)$ choices of the subset $A$, we can sum over these estimates since there are only finitely many of them, since both $\rho_0$ and $k_0$ are fixed constants depending only on $\psi$.
When $p=2$, inequality~\eqref{sufficeest} is a consequence of the chain
\begin{align}\label{shiftl2est}
\Vert T^{y}f\Vert_{L^2(\bbrn)}&\leq\bigg\Vert  \Big(\sum_{k\in \Z}\big|(\psi_{k})^y\ast f \big|^2 \Big)^{1/2} \bigg\Vert_{L^2(\bbrn)}
=\bigg( \sum_{k\in \Z}\big\Vert (\psi_{k})^y\ast (\wt{\psi}_{k}\ast f)\big\Vert_{L^2(\bbrn)}^2\bigg)^{1/2}\\
&\lesssim \bigg( \sum_{k\in\bbz}\big\Vert \wt{\psi}_{k}\ast f \big\Vert_{L^2(\bbrn)}^2\bigg)^{1/2}\nonumber
\sim\Vert f\Vert_{L^2(\bbrn)}.
\end{align}
Moreover, if $l \in \Z$ then \eqref{doconest} yields
\begin{align}\label{psikasttyfx}
\big| \psi_l\ast T^{y}f(x)\big|
\le \sum_{k\in\bbz}\big|\psi_l\ast (\psi_k)^y\ast f(x) \big|   \lesssim_{\sigma} \sum_{k:2^k\sim 2^l}\mathfrak{M}_{\sigma,2^k}\big( (\psi_k)^y\ast f\big)(x)
\end{align}
for $\sigma>n$.
By using \eqref{psikasttyfx}, \eqref{bmomaximaline}, and \eqref{suppdbmoest},
\begin{align}\label{tyfbmoest}
\big\Vert T^{y}f\big\Vert_{BMO(\bbrn)}&\sim \sup_{P\in\mathcal{D}}\bigg(\frac{1}{|P|}\int_P{\sum_{l=-\log_2{\ell(P)}}^{\infty}{\big| \psi_l\ast T^{y}f(x)\big|^2     }}\; dx \bigg)^{1/2}\nonumber\\
&\lesssim_{\sigma}  \sup_{R\in\mathcal{D}}\bigg(\frac{1}{|R|}\int_R{\sum_{k=-\log_2{\ell(R)}}^{\infty}{\big| (\psi_k)^y\ast f(x)\big|^2     }}\; dx \bigg)^{1/2}\nonumber\\
&\lesssim \big( \log{(e+|y|)}\big)^{1/2}  \Vert f\Vert_{BMO(\bbrn)}.
\end{align}
Interpolating \eqref{shiftl2est} and \eqref{tyfbmoest}, for $2<p<\infty$
\begin{align*}
 \big\Vert T^{y}f\big\Vert_{L^p(\bbrn)}\lesssim \big( \log{(e+|y|)}\big)^{1/2-1/p}\Vert f\Vert_{L^p(\bbrn)},
\end{align*} as desired.
For the case $1\le p<2$, we employ duality arguments with the transpose operator 
$$(T^{y})^*f(x)=\sum_{k\in A}\overline{(\psi_{k})^{-y}}\ast f(x).$$
Similarly to \eqref{tyfbmoest}, we have
$$\big\Vert (T^{y})^*f\big\Vert_{BMO(\bbrn)}\lesssim \big( \log{(e+|y|)}\big)^{1/2}\Vert f\Vert_{BMO(\bbrn)},$$
which yields that
\[
 \big\|T^{y}f \big\|_{H^1(\R^n)} \lesssim \big( \log{(e+|y|)}\big)^{1/2}\Vert f\Vert_{H^1(\bbrn)}. 
\]
By interpolating this inequality with \eqref{shiftl2est}, we finally obtain \eqref{sufficeest} for $1\le p<2$. 
This completes the proof.
\end{proof}

\section{Proof of Proposition \ref{T:Hormander}}\label{S:linear}

In this section, we prove Proposition~\ref{T:Hormander}.
Let $\vartheta$ be a radial Schwartz function having the properties that
\begin{equation}\label{psiftpro}
\supp{(\wh{\vartheta})}\subset \{\xi\in\bbrn:1/4\le |\xi|\le 4\} \q \text{ and }\q \wh{\vartheta}(\xi)=1 ~\text{ for }~1/2\le |\xi|\le 2,
\end{equation}
and
$$\sum_{k\in\bbz}\wh{\vartheta_k}(\xi)=2 ~\text{ for }~ \xi\not=0.$$

\subsection{Proof of \eqref{E:opt-est}}
Assume that $\lambda \geq |1/2-1/p|$. 
Let $\Lambda$ be the bilinear form given by
\[
\Lambda(f_1,f_2):=\int_{\R^n} Tf_1(x) f_{2}(x)\,dx,
\]
for $f_1, f_2\in \mathscr{S}_0(\bbrn)$.
Since 
$$(H^1(\R^n))^*=BMO(\R^n) \q \text{ and }\q (L^p(\R^n))^*=L^{p'}(\R^n)~ \text{ for}~ 1<p<\infty,$$
the estimate~\eqref{E:opt-est} will follow if we prove that
\begin{equation}\label{E:form}
\big| \Lambda(f_1,f_2) \big| \leq CD_\lambda(K) \|f_1\|_{X^{p}(\R^n)} \|f_{2}\|_{X^{p'}(\R^n)}.
\end{equation}
We rewrite the form $\Lambda$ as follows:
	\begin{align*} \Lambda (f_1,f_2)
		&= \sum_{k\in \Z} \int_{\R^n} \int_{\R^{n}} \wh {K}(2^{-k}\xi\,)  \wh{f_1}(\xi)  f_{2}(x) e^{2\pi i \langle x,\xi\rangle} \; d\xi dx \\
		&= \sum_{k\in \Z} \int_{\R^{n}} \wh {K}(2^{-k}\xi)   \wh f_1(\xi) \wh{f_{2}}(-\xi)\; d\xi.
	\end{align*}
The previous expression can be further rewritten as
\begin{equation*}
	\sum_{k\in \Z} \int_{\R^{n}} \wh {K}(2^{-k}\xi\,)  \wh f_1(\xi)\wh{\vartheta}(2^{-k}\xi)  \wh{f_{2}}(-\xi) \wh{\vartheta}(-2^{-k}\xi)\; d\xi,
\end{equation*}
as the Fourier transform of $\vartheta$ coincides with $1$ in $\supp{(\wh{K}})$.
Using Fourier inversion, this is equal to
\begin{align*}& \sum_{k\in \Z} \int_{\R^n} \int_{\R^{n}} K_k (y)   f_1\ast \vartheta_{k}(x-y)  f_{2}\ast \vartheta_{k}(x) \;dydx \\
= &\sum_{k\in \Z} \int_{\R^n} \int_{\R^{n}}  {K}(y)  f_1\ast \vartheta_{k}(x-2^{-k}y)  f_{2}\ast \vartheta_{k}(x) \;dydx \\
= &\sum_{k\in \Z} \int_{\R^n} {K}(y)  \Big( \int_{\R^{n}}    f_1\ast (\vartheta_{k})^{y }(x) f_2\ast \vartheta_{k}(x) \; dx\Big)\;dy.
\end{align*}
Therefore,
\begin{equation*}
\big| \Lambda(f_1,f_2)\big|\le \int_{\bbrn}|K(y)|\Big(\int_{\bbrn}\sum_{k\in\bbz} \mathfrak{M}_{\sigma,2^k}\big( (\vartheta_k)^y\ast f_1\big)(x)\mathfrak{M}_{\sigma,2^k}\big(\vartheta_k\ast f_2\big)(x) \; dx\Big)\; dy
\end{equation*} for $\sigma>n$.
Now we fix $y\in\bbrn$ in order to estimate the innermost integral. By applying Corollary \ref{C:ShiftedSquare}, for each dyadic cube $Q\in\mathcal{D}$ we choose proper measurable subsets $S_Q^1$ and $S_Q^2$ of $Q$ such that
\begin{equation}\label{ssq12ine}
|S_Q^1|, |S_Q^2|>\frac{3}{4}|Q|,
\end{equation}
\begin{equation}\label{sestf2}
\bigg\Vert \Big\{ \sum_{Q\in\mathcal{D}_k}\Big( \inf_{z\in Q}\mathfrak{M}_{\sigma,2^k}{\big((\vartheta_k)^y\ast f_1\big)}(z)\Big)\chi_{S_Q^1}\Big\}_{k\in\bbz}\bigg\Vert_{L^p(\ell^2)}\lesssim \big(\log{(e+|y|)} \big)^{|1/p-1/2|}\Vert f_1 \Vert_{X^{p}(\bbrn)},
\end{equation}
and
\begin{equation}\label{sestf3}
\bigg\Vert \Big\{ \sum_{Q\in\mathcal{D}_k}\Big( \inf_{z\in Q}\mathfrak{M}_{\sigma,2^k}{\big(\vartheta_k\ast f_2\big)}(z)\Big)\chi_{S_Q^2}\Big\}_{k\in\bbz}\bigg\Vert_{L^{p'}(\ell^2)}\lesssim \Vert f_2 \Vert_{X^{p'}(\bbrn)}.
\end{equation} 
We observe that \eqref{ssq12ine} implies
$$|S_Q^1\cap S_Q^2|>\frac{1}{2}|Q|.$$
Then we estimate, using also~\eqref{infmax2},
\begin{align*}
&\int_{\bbrn}\sum_{k\in\bbz} \mathfrak{M}_{\sigma,2^k}\big( (\vartheta_k)^y\ast f_1\big)(x)\mathfrak{M}_{\sigma,2^k}\big(\vartheta_k\ast f_2\big)(x) \; dx\\
&= \int_{\bbrn}\sum_{k\in\bbz}\sum_{Q\in\mathcal{D}_k} \mathfrak{M}_{\sigma,2^k}\big( (\vartheta_k)^y\ast f_1\big)(x)\mathfrak{M}_{\sigma,2^k}\big(\vartheta_k\ast f_2\big)(x)\chi_Q(x) \; dx\\
&\sim \sum_{k\in\bbz}\sum_{Q\in\mathcal{D}_k}\Big(\inf_{z\in Q}\mathfrak{M}_{\sigma,2^k}\big( (\vartheta_k)^y\ast f_1\big)(z)\Big)\Big(\inf_{z\in Q}\mathfrak{M}_{\sigma,2^k}\big(\vartheta_k\ast f_2\big)(z)  \Big)|Q|\\
&\sim \sum_{k\in\bbz}\sum_{Q\in\mathcal{D}_k}\Big(\inf_{z\in Q}\mathfrak{M}_{\sigma,2^k}\big( (\vartheta_k)^y\ast f_1\big)(z)\Big)\Big(\inf_{z\in Q}\mathfrak{M}_{\sigma,2^k}\big(\vartheta_k\ast f_2\big)(z)  \Big)|S_Q^1\cap S_Q^2|\\
&\le \int_{\bbrn}\sum_{k\in\mathbb{Z}}\bigg( \sum_{Q\in\mathcal{D}_k}\Big(\inf_{z\in Q}\mathfrak{M}_{\sigma,2^k}\big( (\vartheta_k)^y\ast f_1\big)(z)\Big)\chi_{S_Q^1}(x)\bigg)\\
&\qq\qq\qq\times \bigg( \sum_{Q\in\mathcal{D}_k}\Big(\inf_{z\in Q}\mathfrak{M}_{\sigma,2^k}\big( \vartheta_k\ast f_2\big)(z)\Big)\chi_{S_Q^2}(x)\bigg)\; dx.
\end{align*}
Now H\"older's inequality, \eqref{sestf2} and \eqref{sestf3} prove that the last expression is bounded by
\begin{align*}
&\Big\Vert \Big\{  \sum_{Q\in\mathcal{D}_k}\Big(\inf_{z\in Q}\mathfrak{M}_{\sigma,2^k}\big( (\vartheta_k)^y\ast f_1\big)(z)\Big)\chi_{S_Q^1}\Big\}_{k\in\bbz}\Big\Vert_{L^{p}(\ell^2)}\\
&\qq\times \Big\Vert \Big\{ \sum_{Q\in\mathcal{D}_k}\Big(\inf_{z\in Q}\mathfrak{M}_{\sigma,2^k}\big( \vartheta_k\ast f_2\big)(z)\Big)\chi_{S_Q^2} \Big\}_{k\in\bbz}\Big\Vert_{L^{p'}(\ell^2)}\\
&\lesssim \big(\log{(e+|y|)} \big)^{|1/p-1/2|}\Vert f_1 \Vert_{X^{p}(\bbrn)}\Vert f_2 \Vert_{X^{p'}(\bbrn)}.
\end{align*}
This yields~\eqref{E:form}.

\subsection{Optimality of the condition $\lambda \geq |1/2-1/p|$}\label{SS:OptimalityLinear}

Assume that inequality~\eqref{E:opt-est} holds for some $C>0$, for all kernels $K$ and $f\in \mathscr{S}_0(\bbrn)$. We will show that then necessarily $\lambda \geq |1/2-1/p|$. For each $k\in\bbz$, let $\zeta_k:=10k$.

We first assume that $2\le p\le \infty$. 
 Let  $\eta$ and $\beta$ be radial Schwartz functions such that
 \begin{equation}\label{cutoffcon1}
 \supp{(\wh{\eta})}\subset \Big\{\xi\in\bbrn:|\xi|\le \frac{1}{100}\Big\}, \q \supp{(\wh{\beta})}\subset \Big\{\xi\in\bbrn: \frac{10}{11}\le |\xi|\le \frac{11}{10} \Big\},
 \end{equation}
 and
 \begin{equation}\label{cutoffcon2}
 \wh{\beta}(\xi)=1 \q \text{ for }~ \frac{20}{21}\le |\xi|\le \frac{21}{20}.
 \end{equation}
  We fix $N\in \N$ and set
\[
K(y):=\beta(y-2^{\zeta_N} e_1)
\]
and
\[
f(x):=\sum_{k=1}^N \eta(x+2^{\zeta_N-\zeta_k}e_1) e^{2\pi i \langle x, 2^{\zeta_k} e_1\rangle}.
\]
Then
\[
\wh{K}(\xi)=\wh{\beta}(\xi) e^{-2\pi i \langle 2^{\zeta_N} e_1,\xi\rangle} \q \text{and} \q \wh{f}(\xi)=\sum_{k=1}^N \wh{\eta}(\xi-2^{\zeta_k} e_1) e^{2\pi i \langle  2^{\zeta_N-\zeta_k} e_1,\xi\rangle}.
\]
Given $l\in \mathbb Z$, the assumptions \eqref{cutoffcon1} and \eqref{cutoffcon2} yield that
\begin{equation}\label{tfcombineest}
\wh{\beta}(\xi/2^l)\wh{\eta}(\xi-2^{\zeta_k}e_1)=\begin{cases} 0 & \text{if } l\not= {\zeta_k},\\
\wh{\eta}(\xi-2^{\zeta_k}e_1) & \text{if } l={\zeta_k},\end{cases}
\end{equation} 
and thus
\begin{equation*}
\wh{Tf}(\xi)=\sum_{k=1}^{N}\wh{\eta}(\xi-2^{\zeta_k}e_1),
\end{equation*}
and equivalently,
\begin{equation}\label{tfxre}
Tf(x)=\sum_{k=1}^N \eta(x) e^{2\pi i \langle x, 2^{\zeta_k} e_1\rangle}.
\end{equation}

Now we claim that
\begin{equation}\label{coficon}
\Vert f\Vert_{X^p(\bbrn)}\lesssim N^{1/p}\q \text{ and } \q \Vert Tf\Vert_{X^p(\bbrn)}\gtrsim N^{1/2}.
\end{equation}
Since
\begin{equation}\label{dlamkjlam}
D_\lambda(K) \sim N^\lambda,
\end{equation}
 the hypothetical inequality~\eqref{E:opt-est}, together with \eqref{coficon}, finally yields
\[
N^{{1}/{2}} \lesssim N^{{1}/{p}+\lambda}.
\]
Sending $N$ to infinity, we obtain $\lambda \geq 1/2-1/p$, as desired.
Therefore, it remains to verify \eqref{coficon}.
To this end,
we observe that
\begin{equation}\label{psilastf1}
\vartheta_l\ast f=\begin{cases}\vartheta_l\ast \big(\eta(\cdot+2^{\zeta_N-\zeta_k}e_1)e^{2\pi i\langle \cdot,2^{\zeta_k}e_1\rangle} \big) & \text{if }\zeta_k-2\le l\le \zeta_k+2,~ 1\le k\le N,\\
0 & \text{otherwise},
\end{cases}
\end{equation}
where $\vartheta$ is the radial Schwartz function defined at the beginning of this section.
Therefore,  for each $1\le k\le N$ and $\zeta_k-2\le l\le \zeta_k+2$,
\begin{equation}\label{psilastf2}
\big| \vartheta_l\ast f(x)\big|\le |\vartheta_l|\ast \big| \eta(\cdot+2^{\zeta_N-\zeta_k}e_1)\big|(x)\lesssim \frac{1}{(1+|x+2^{\zeta_N-\zeta_k}e_1|)^n}.
\end{equation}

When $2\le p<\infty$, using the Littlewood-Paley theorem \eqref{hardycha}, we have
\begin{align}\label{estinputf}
\Vert f\Vert_{L^p(\bbrn)}&\sim\Big\Vert \Big(\sum_{l\in\bbz}\big| \vartheta_l\ast f\big|^2 \Big)^{1/2}\Big\Vert_{L^p(\bbrn)}\nonumber\\
&\lesssim \bigg(\int_{\bbrn}    \Big(  \sum_{k=1}^{N}\frac{1}{(1+|x+2^{\zeta_N-\zeta_k}e_1|)^{2n}}      \Big)^{p/2}        \, dx \bigg)^{1/p}\nonumber\\
&\lesssim  \bigg(\int_{\bbrn}  \sum_{k=1}^{N}  \frac{1}{(1+|x+2^{\zeta_N-\zeta_k}e_1| )^{2n}}    \,dx \bigg)^{1/p}\nonumber\\
&\lesssim N^{1/p},
\end{align}
where the penultimate estimate follows from the fact that
\begin{equation}\label{sumk1jfra}
\sum_{k=1}^{N}  \frac{1}{(1+|x+2^{\zeta_N-\zeta_k}e_1| )^{2n}}\lesssim 1 \q \text{ uniformly in }~ x\in\bbrn.
\end{equation}

On the other hand, using \eqref{hardycha} again,
\begin{align*}
\Vert Tf\Vert_{L^p(\bbrn)}\sim\Big\Vert \Big(\sum_{l\in\bbz}\big| \vartheta_l\ast Tf\big|^2 \Big)^{1/2}\Big\Vert_{L^p(\bbrn)}\ge\Big\Vert \Big(\sum_{k=1}^{N}\big| \vartheta_{\zeta_k}\ast Tf\big|^2 \Big)^{1/2}\Big\Vert_{L^p(\bbrn)}
\end{align*}
and, in view of \eqref{psiftpro} and \eqref{tfxre}, we see
$$\big|\vartheta_{\zeta_k}\ast Tf(x)\big|=\big| \eta(x)\big|,$$
which finally yields that
$$\Vert Tf\Vert_{L^p(\bbrn)}\gtrsim N^{1/2}\Vert \eta\Vert_{L^p(\bbrn)}\sim N^{1/2}.$$
This proves \eqref{coficon} for $2\le p<\infty$.

When $p=\infty$, using \eqref{bmocha}, \eqref{psilastf1}, \eqref{psilastf2}, and \eqref{sumk1jfra} we have
\begin{align*}
\Vert f\Vert_{BMO(\bbrn)}&\sim \sup_{P\in\mathcal{D}}\bigg(\frac{1}{|P|}\int_P   \sum_{l=-\log_2{\ell(P)}}^{\infty}\big| \vartheta_l\ast f(x)\big|^2    \, dx \bigg)^{1/2}\\
&\lesssim \sup_{P\in\mathcal{D}}\Big( \frac{1}{|P|}\int_P \sum_{k=1}^N \frac{1}{(1+|x+2^{\zeta_N-\zeta_k}e_1|)^{2n}}    \; dx\Big)^{1/2}\\
&\lesssim \sup_{P\in\mathcal{D}}\Big(\frac{1}{|P|}\int_P 1\; dx \Big)^{1/2}\le 1.
\end{align*}
In addition,
\begin{align*}
\Vert Tf\Vert_{BMO(\bbrn)}&\gtrsim \bigg(\int_{[0,1]^n}\sum_{l=0}^{\infty}\big| \vartheta_l\ast Tf(x)\big|^2 \, dx\bigg)^{1/2}\ge \bigg(\int_{[0,1]^n}\sum_{k=1}^{N}\big| \vartheta_{\zeta_k}\ast Tf(x)\big|^2 \, dx\bigg)^{1/2}\\
&=\Big(\int_{[0,1]^n} \sum_{k=1}^{N} |\eta(x)|^2   \; dx \Big)^{1/2}=N^{1/2}\Vert \eta\Vert_{L^2([0,1]^n)}\sim N^{1/2},
\end{align*}
which completes the proof of \eqref{coficon}.

Let us next assume that $1\le p<2$. Let $T^{*}$ denote the transpose operator of $T$. 
This operator is of the same type as $T$ but associated with the kernel $\wt{K}(y)=K(-y)$. As $D_\lambda(K)=D_\lambda(\wt{K})$, the inequality~\eqref{E:opt-est} is equivalent to 
\[
\|T^{*}g\|_{X^{p'}(\R^n)} \leq C D_\lambda(\wt{K}) \|g\|_{X^{p'}(\R^n)}. 
\]
Since $2<p'\le \infty$, the above inequality fails to hold for $\lambda<1/2-1/p'=1/p-1/2$ by the result of the previous paragraph.

\section{Proof of Theorem~\ref{T:HormanderB}}\label{S:bilinear}

In this section, we prove \Cref{T:HormanderB}.

\subsection{Proof of \eqref{E:BoundB}}\label{ss:sufficiency-bilinear}
Let $\Lambda$ be the trilinear form given by
\[
\Lambda(f_1,f_2,f_3):=\int_{\R^n} \mathcal{B}\big(f_1,f_2\big)(x) f_{3}(x)\,dx,
\]
 for $f_1,f_2,f_3\in\mathscr{S}_0(\R^n)$.
We observe that \eqref{E:BoundB} will follow if we prove that
\begin{equation}\label{E:formB}
\big| \Lambda(f_1,f_2,f_3) \big| \lesssim_{p_1,p_2,n} D_\lambda(K) \|f_1\|_{H^{p_1}(\bbrn)} \|f_2\|_{H^{p_2}(\bbrn)} \|f_{3}\|_{H^{p'}(\bbrn)}
\end{equation}
for $\lambda\geq\max\Big\{\frac1{p_1}, \frac{1}{p_2}, \frac{1}{p'}\Big\}$. 

The transpose operators $\mathcal{B}^{\ast_1}$ and $\mathcal{B}^{\ast_2}$ are defined so that 
\[
\Lambda(f_1,f_2,f_3)= \int_{\R^n} \mathcal{B}^{\ast_1}\big(f_3,f_2\big)(x) f_{1}(x)\,dx = \int_{\R^n} \mathcal{B}^{\ast_2}\big(f_1,f_3\big)(x) f_{2}(x)\,dx,
\]
and they are operators of the same form as $\mathcal{B}$, but with initial kernels
\[
K^{\ast_1}(y_1,y_2) = K(-y_1, -y_1 + y_2) \quad \text{and} \quad K^{\ast_2}(y_1,y_2) = K(y_1-y_2, -y_2),
\]
respectively. Note that both $K^{\ast_1}$ and $K^{\ast_2}$ satisfy similar properties as $K$, in the sense that
\begin{equation}\label{E:kernel-equivalence}
D_\lambda(K^{\ast_1}) \sim_{n} D_\lambda(K^{\ast_2}) \sim_{n} D_\lambda(K) \end{equation}
and
\[
\wh{K^{\ast_1}} \text{ and } \wh{K^{\ast_2}}  \text{ are supported on } |(\xi_1,\xi_2)|\sim_{n} 1.
\] 

Using Fourier inversion, we rewrite the form $\Lambda$ as
	\begin{align}\label{E:LambdaRewrite} &\Lambda(f_1,f_2,f_3)
		= \sum_{k\in \Z} \int_{\R^{2n}} \wh {K}(2^{-k}\xi_1,2^{-k}\xi_2)   \wh f_1(\xi_1) \wh{f_2}(\xi_2)\wh{f_{3}}(-\xi_1-\xi_2) \;d\xxi,
	\end{align}
	where $d\xxi=d\xi_1 d\xi_2$.

We next decompose $\wh{K}$ into pieces supported in certain cones in the frequency space. Since $\wh{K}$ is supported in the unit annulus, for each $(\xi_1,\xi_2)$ in the support of $\wh{K}$, at least two of the three terms $|\xi_1|$, $|\xi_2|$ and $|-\xi_1 - \xi_2|$ are bounded away from zero. By a smooth partition of unity, it therefore suffices to consider pieces of the form
\begin{equation}\label{E:decomp}
\wh{K}(\xi_1,\xi_2) \wh{\Phi_1}(\xi_1) \wh{\Phi_2}(\xi_2) \wh{\Phi_3}(-\xi_1-\xi_2),
\end{equation} 
where two of the functions $\Phi_i$ coincide with $\psi$ and the remaining one is equal to $\phi$,
with  $\phi$ and $\psi$ being the functions introduced at the beginning of Section \ref{S:preliminary}.
We rewrite those terms in~\eqref{E:decomp} for which $\Phi_3=\psi$ more symmetrically as 
\begin{equation}\label{E:transpose-kernel}
\wh{K^{\ast_1}}(-\xi_1-\xi_2,\xi_2) \wh{\psi}(-\xi_1-\xi_2) \wh{\psi}(\xi_2) \wh{\phi}(\xi_1) \quad \text{ and } \quad \wh{K^{\ast_2}}(\xi_1,-\xi_1-\xi_2) \wh{\psi}(\xi_1) \wh{\psi}(-\xi_1-\xi_2) \wh{\phi}(\xi_2).
\end{equation}

Decomposing the form~\eqref{E:LambdaRewrite} into three pieces according to the decomposition of $\wh{K}$ described above, and performing the change of variables in the pieces involving~\eqref{E:transpose-kernel}, it then suffices to estimate the terms
\begin{align*}
	\Lambda_0(f_1,f_2,f_3) &:= \sum_{k\in \Z} \int_{\R^{2n}} \wh {K}(2^{-k}\xi_1,2^{-k}\xi_2)  \wh{\psi_k\ast f_1}(\xi_1)  \wh{\psi_k\ast f_{2}}(\xi_2) \wh{\phi_k\ast f_{3}}(-\xi_1-\xi_2)\; d\xxi,\\
    \Lambda_1(f_1,f_2,f_3) &:= \sum_{k\in \Z} \int_{\R^{2n}} \wh {K^{\ast_1}}(2^{-k}\xi_1,2^{-k}\xi_2)  \wh{\psi_k\ast f_3}(\xi_1)  \wh{\psi_k\ast f_{2}}(\xi_2) \wh{\phi_k\ast f_{1}}(-\xi_1-\xi_2) \; d\xxi,\\
    \Lambda_2(f_1,f_2,f_3) &:= \sum_{k\in \Z} \int_{\R^{2n}} \wh {K^{\ast_2}}(2^{-k}\xi_1,2^{-k}\xi_2)\,)  \wh{\psi_k\ast f_1}(\xi_1)  \wh{\psi_k\ast f_{3}}(\xi_2) \wh{\phi_k\ast f_{2}}(-\xi_1-\xi_2)\; d\xxi.
\end{align*}

We focus on estimating $\Lambda_0(f_1,f_2,f_3)$, as symmetric arguments can be applied to the other two cases. 
Using Fourier inversion again, we can write
\begin{equation*}
\Lambda_0(f_1,f_2,f_3) = \int_{\R^{2n}}K(\yyy)\Big( \sum_{k\in \Z}  \int_{\R^{n}}    (\psi_k)^{y_1}\ast f_1(x) (\psi_k)^{y_2}\ast f_2(x) \phi_k\ast f_3(x) \; dx \Big) \;d\yyy.
\end{equation*}
Now we claim that for each $\yyy\in\bbr^{2n}$,
\begin{align}\label{la0red}
&\bigg| \sum_{k\in \Z}  \int_{\R^{n}}    (\psi_k)^{y_1}\ast f_1(x) (\psi_k^{y_2})\ast f_2(x) \phi_k\ast f_3(x) \; dx\bigg|\nonumber\\
&\lesssim \log(e+|\yyy|)^{|1/2-1/{p_1}|+|1/2-1/{p_2}|} \Vert f_1\Vert_{H^{p_1}(\bbrn)}\Vert f_2\Vert_{H^{p_2(\bbrn)}}\Vert f_3\Vert_{H^{p'}(\bbrn)}.
\end{align}
We fix $\sigma>n$ and use the notation
\begin{equation*}
\mathfrak{V}_{\sigma,k}f:=\mathfrak{M}_{\sigma,2^k}\big(\phi_k\ast f\big)\q \text{ and }\q \mathfrak{U}_{\sigma,k}^{y}f:=\mathfrak{M}_{\sigma,2^k}\big((\psi_k)^y\ast f\big).
\end{equation*}
To verify \eqref{la0red}, we first apply \Cref{C:ShiftedSquare} to choose two sequences $\{S_Q^1\}_{Q\in\mathcal{D}}$ and $\{S_Q^2\}_{Q\in\mathcal{D}}$ such that for each $l=1,2$
$$S_Q^l\subset Q, \q |S_Q^l|>\frac{3}{4}|Q|,$$
and 
\begin{equation}\label{frakmaxestsql}
\bigg\Vert \Big\{ \sum_{Q\in\mathcal{D}_k}\Big(\inf_{z\in Q} \mathfrak{U}_{\sigma,k}^{y}f_l(z)   \Big)\chi_{S_Q^l}\Big\}_{k\in\bbz}\bigg\Vert_{L^{p_l}(\ell^2)}\lesssim \big(\log{(e+|y|)}\big)^{|1/p_l-1/2|}\Vert f_l\Vert_{H^{p_l}(\bbrn)}
\end{equation}
where the constants in the inequalities are independent of $f_1$ and $f_2$ while the subsets $S_Q^1$ and $S_Q^2$ depend on $f_1$ and $f_2$, respectively.
We note that 
\begin{equation}\label{sq1sq2estq}
|S_Q^1\cap S_Q^2 |>\frac{1}{2}|Q|.
\end{equation}
Now the left-hand side of \eqref{la0red} is bounded by
\begin{align*}
&\sum_{k\in\bbz}\sum_{Q\in\mathcal{D}_k}\int_{Q}\mathfrak{U}_{\sigma,k}^{y_1}f_1(x) \mathfrak{U}_{\sigma,k}^{y_2}f_2(x)\mathfrak{V}_{\sigma,k}f_3(x)\; dx\\
&\sim \sum_{k\in\bbz}\sum_{Q\in\mathcal{D}_k}\Big(\inf_{z\in Q} \mathfrak{U}_{\sigma,k}^{y_1}f_1(z)\Big) \Big(\inf_{z\in Q}\mathfrak{U}_{\sigma,k}^{y_2}f_2(z)\Big) \Big(\inf_{z\in Q}\mathfrak{V}_{\sigma,k}f_3(z)\Big) |Q|\\
&\sim \sum_{k\in\bbz}\sum_{Q\in\mathcal{D}_k}\Big(\inf_{z\in Q} \mathfrak{U}_{\sigma,k}^{y_1}f_1(z)\Big) \Big(\inf_{z\in Q}\mathfrak{U}_{\sigma,k}^{y_2}f_2(z)\Big) \Big(\inf_{z\in Q}\mathfrak{V}_{\sigma,k}f_3(z)\Big)|S_Q^1\cap S_Q^2|\\
&\le \sum_{k\in\bbz}\int_{\bbrn}\mathfrak{V}_{\sigma,k}f_3(x)\sum_{Q\in\mathcal{D}_k}\Big(\inf_{z\in Q} \mathfrak{U}_{\sigma,k}^{y_1}f_1(z)\Big)\Big(\inf_{z\in Q} \mathfrak{U}_{\sigma,k}^{y_2}f_2(z)\Big)\chi_{S_Q^1\cap S_Q^2}(x) \; dx,
\end{align*}
where \eqref{infmax2} and \eqref{sq1sq2estq} are applied.
By using H\"older's inequality, \eqref{peetremax} and \eqref{frakmaxestsql}, the last expression is estimated by
\begin{align*}
&\int_{\bbrn} \big\Vert\big\{\mathfrak{V}_{\sigma,k}f_3(x) \big\}_{k\in\bbz}\big\Vert_{\ell^{\infty}}  \prod_{l=1}^{2} \Big\Vert \Big\{ \sum_{Q\in\mathcal{D}_k}\Big(\inf_{z\in Q} \mathfrak{U}_{\sigma,k}^{y_l}f_l(z)\Big)\chi_{S_Q^l}(x)\Big\}_{k\in\bbz}\Big\Vert_{\ell^2}\; dx\\
&\le \Big\Vert\big\{\mathfrak{V}_{\sigma,k}f_3 \big\}_{k\in\bbz}\Big\Vert_{L^{p'}(\ell^{\infty})} \prod_{l=1}^{2}\bigg\Vert \Big\{ \sum_{Q\in\mathcal{D}_k}\Big( \inf_{z\in Q} \mathfrak{U}_{\sigma,k}^{y_l}f_l(z)\Big)\chi_{S_Q^l}\Big\}_{k\in\bbz}\bigg\Vert_{L^{p_l}(\ell^2)}\\
&\lesssim        \log(e+|y_1|)^{|1/p_1-1/2|}  \log(e+|y_2|)^{|1/p_2-1/2|}    \Vert f_1\Vert_{H^{p_1}(\bbrn)} \Vert f_2\Vert_{H^{p_2}(\bbrn)}\Vert f_3\Vert_{H^{p'}(\bbrn)},
\end{align*}
which completes the proof of \eqref{la0red}.
Finally, \eqref{la0red} implies that for $1\le p_1,p_2,p'\le \infty$
\begin{equation*}
|\Lambda_0(f_1,f_2,f_3)| \leq CD_\lambda(K) \|f_1\|_{H^{p_1}(\bbrn)} \|f_2\|_{H^{p_2}(\bbrn)} \|f_{3}\|_{H^{p'}(\bbrn)}
\end{equation*}
if $\lambda\geq |1/2-1/{p_1}|+|1/2-1/{p_2}|$.

 The same argument can be used to estimate $\Lambda_1(f_1,f_2,f_3)$ and $\Lambda_2(f_1,f_2,f_3)$. Thanks to~\eqref{E:kernel-equivalence}, these forms are bounded by the right-hand side of~\eqref{E:formB} if $\lambda\geq |1/2-1/{p'}|+|1/2-1/{p_2}|$ and $\lambda\geq |1/2-1/{p_1}|+|1/2-1/{p'}|$, respectively. Thus,~\eqref{E:formB} holds if 
\begin{align*}
\lambda &\geq \max\bigg\{\Big|\frac12-\frac1{p_1}\Big|+\Big|\frac12-\frac1{p_2}\Big|,\,\Big|\frac12-\frac1{p'}\Big|+\Big|\frac12-\frac1{p_2}\Big|,\, \Big|\frac12-\frac1{p_1}\Big|+\Big|\frac12-\frac1{p'}\Big| \bigg\}\\
&= \max\bigg\{\frac1{p_1}, \frac{1}{p_2}, \frac{1}{p'}\bigg\}.
\end{align*}

\subsection{Optimality of the condition $\lambda\ge \max\{\frac{1}{p_1},\frac{1}{p_2},\frac{1}{p'}\}$}

We split the range of exponents $p_1$, $p_2$ into three regions depending on which of the numbers $\frac{1}{p_1}$, $\frac{1}{p_2}$ and $\frac{1}{p'}$ is the largest, see~\Cref{fig:B2} below. 
Due to the symmetry between the operator $\mathcal B$ and its transpose operators $\mathcal{B}^{\ast_1}$ and $\mathcal{B}^{\ast_2}$, described at the beginning of Subsection~\ref{ss:sufficiency-bilinear}, a counterexample in one of the three regions implies counterexamples in the other two. Thus, the proof will be complete if we show that whenever \begin{equation}\label{conditionpp12}
\max\Big\{\frac{1}{p_1},\frac{1}{p_2},\frac{1}{p'}\Big\}=\frac{1}{p'}
\end{equation} and if equation~\eqref{E:BoundB} holds for all kernels $K$, then necessarily $\lambda \geq \frac{1}{p'}$. 
We note that the condition \eqref{conditionpp12} implies $p_1,p_2\ge 2$.

\begin{figure}[h]
\centering
\scalebox{1.5}{\includegraphics[page=1]{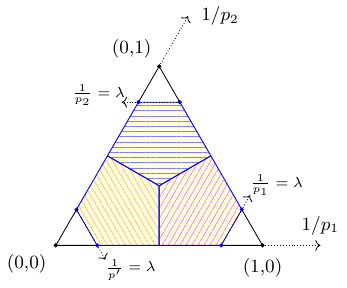}}
\caption{Range of $p_1,p_2$ for $H^{p_1}\times H^{p_2}\to Y^p$ boundedness of $\mathcal{B}$ for a given value of $\lambda$, split into three regions depending on which of the numbers $\frac{1}{p_1}$, $\frac{1}{p_2}$ and $\frac{1}{p'}$ is the largest.}
\label{fig:B2}
\end{figure}

We apply similar arguments as in Section \ref{SS:OptimalityLinear}.
Let  $\eta$ and $\beta$ be nontrivial radial Schwartz functions satisfying~\eqref{cutoffcon1} and~\eqref{cutoffcon2}.
For each $k\in\bbz$, let $\zeta_k:=10k$.
We fix $N\in\bbn$ and set
\begin{align*}
f(x)&:=\sum_{k=1}^{N}\eta(x+2^{\zeta_N-\zeta_k}e_1)e^{2\pi i\langle x,2^{\zeta_k} e_1\rangle},\\
 g(x)&:=\sum_{k=1}^{N}\eta(x+2^{\zeta_N-\zeta_k}e_1)e^{-2\pi i\langle x,2^{\zeta_k} e_1\rangle},
\end{align*}
and
$$K(y_1,y_2):=\beta(y_1-2^{\zeta_N} e_1)\beta(y_2-2^{\zeta_N} e_1).$$
Then we see that
\begin{align*}
\wh{K}(\xi_1,\xi_2)&= \wh{\beta} (\xi_1) \wh{\beta}(\xi_2) e^{-2\pi i\langle 2^{\zeta_N}e_1,\xi_1+\xi_2\rangle},\\
\wh{f}(\xi)&=\sum_{k=1}^{N}\wh{\eta}(\xi-2^{\zeta_k} e_1)e^{2\pi i\langle 2^{\zeta_N-\zeta_k}e_1,\xi\rangle},\\
\wh{g}(\xi)&=\sum_{k=1}^{N}\wh{\eta}(\xi+2^{\zeta_k} e_1)e^{2\pi i\langle 2^{\zeta_N-\zeta_k}e_1,\xi\rangle}.
\end{align*}
If $p_1<\infty$, then \eqref{estinputf} shows that
$$\Vert f\Vert_{L^{p_1}(\bbrn)}\lesssim N^{1/p_1}.$$
Moreover,
\[
\|f\|_{L^\infty(\R^n)} \le \sup_{x\in \R^n} \sum_{k=1}^{N} |\eta(x+2^{\zeta_N-\zeta_k}e_1)|\lesssim \sup_{x\in\bbrn}\sum_{k=1}^N\frac{1}{1+|x+2^{\zeta_N-\zeta_k}e_1|} \lesssim 1.
\]
Similarly, we have 
\[\Vert g\Vert_{L^{p_2}(\bbrn)}\lesssim N^{1/p_2}.\]
On the other hand, using \eqref{tfcombineest}, we have
\begin{align*}
\mathcal{B}(f,g)(x)&= \sum_{k=1}^ N\int_{\R^{2n}}\wh{\eta}(\xi_1-2^{\zeta_k}e_1)\wh{\eta}(\xi_2+2^{\zeta_k}e_1)e^{2\pi i\langle x,\xi_1+\xi_2\rangle}\; d\xi_1 d\xi_2\\
&=\sum_{k=1}^N\big( \eta(x)\big)^2=N \big( \eta(x)\big)^2,
\end{align*}
and thus
\[\big\Vert \mathcal{B}(f,g)\big\Vert_{Y^p(\bbrn)}=N\big\Vert \eta^2\big\Vert_{Y^p(\bbrn)}\sim N.\]
Since $$D_\lambda(K) \sim N^{\lambda},$$
we have
\[\frac{\|\mathcal{B}(f,g)\|_{Y^p(\R^n)}}{D_\lambda(K)\|f\|_{L^{p_1}(\R^n)} \|g\|_{L^{p_2}(\R^n)} }\gtrsim N^{1-1/p - \lambda} = N^{\frac1{p'} - \lambda}, \]
and thus $\lambda \geq \frac1{p'}$ is necessary in order for the condition \eqref{E:BoundB} to hold with a uniform constant $C$.

\section{Proof of Theorem \ref{T:bilinearresult}}\label{proof1}

In this section, we prove Theorem~\ref{T:bilinearresult}. A crucial ingredient of the proof is Theorem~\ref{T:HormanderB}, which is combined with results from the articles~\cite{Gr_He_Ho2018, He_Park2022} and with the classical Coifman-Meyer multiplier theorem~\cite{Co_Me1978} (see also~\cite{Gr_To2002, Ke_St1999}).
The proof is divided into three steps, which are performed in the three subsections below. 

\subsection{Decomposition of the kernel}

Following the approach of Duoandikoetxea and Rubio de Francia \cite{Du_Ru1986}, we reduce the proof of Theorem~\ref{T:bilinearresult} to estimates for operators with smooth kernels. 
Let $\Omega$ be as in the statement of Theorem~\ref{T:bilinearresult}, and 
let $K$ be the kernel of the bilinear operator $T_{\Omega}$ in \eqref{E:rough}, namely
\[
K(\yyy):=\frac{\Omega(\theta)}{|\yyy|^{2n}}, \quad \yyy \neq 0,
\]
where we recall $\theta=\yyy/|\yyy|\in\mathbb{S}^{2n-1}$.
We choose a radial Schwartz function $\Gamma$ on $\bbr^{2n}$ whose Fourier transform $\wh{\Gamma}$ is supported in the annulus $\{\xxxi\in \bbr^{2n}: 1/2\le |\xxxi|\le 2\}$ and satisfies
$\sum_{k\in\bbz}\wh{\Gamma_k}(\xxxi)=1$ for $\xxxi\not= \000$.
For each $k\in\bbz$  we set
 $$ K_{k}(\yyy):=\wh{\Gamma_{-k}}(\yyy)K(\yyy), \quad \yyy\in \R^{2n}$$
 and observe that $K_{k}(\yyy)=2^{2k n} K_0(2^{k} \yyy) = (K_0)_k(\yyy)$.
 Note the slight abuse of notation, namely that $K_k$ is the $L^1$ dilation of $K_0$ instead of $K$.
For each $j\in\bbz$ we define
 \begin{equation}\label{kernelcharacter}
K_{k}^{j}(\yyy):=\Gamma_{j+k}\ast K_{k}(\yyy)=2^{2k n}[\Gamma_j \ast K_{0}](2^{k} \yyy)=2^{2k n}K^{j}_0(2^{k}\yyy) =(K^j_0)_k(\yyy)
\end{equation}
and 
$$K^{j}(\yyy):=\sum_{k \in \bbz}{K^{j}_{k}(\yyy)}.$$
We note that
\begin{equation}\label{whkjxiex}
\wh{K^j}(\xxxi)=\sum_{k\in\bbz}\wh{K^j_0}(2^{-k}\xxxi)=\sum_{k\in\bbz}\wh{K_0}(2^{-k}\xxxi)\wh{\Gamma_{j+k}}(\xxxi).
\end{equation}
The bilinear operator associated with the kernel $K^{j}$ is defined by
 \begin{equation}\label{tjdefinition}
 T^{j}\big(f_1,f_2\big)(x):=\int_{\R^{2n}}{K^{j}(\yyy) f_{1}(x-y_{1})  f_{2}(x-y_{2})} \; d\yyy, \q x\in\bbrn
 \end{equation} 
and the operator $T_{\Omega}$ can be decomposed as
 $$T_{\Omega}=\sum_{j\in\bbz}T^j.$$ 
Thus
 \begin{align}\label{E:series}
 \Vert T\Vert_{L^{p_1}\times L^{p_2}\to L^p}&\lesssim \sum_{j\in\bbz}\Vert T^{j}\Vert_{L^{p_1}\times L^{p_2}\to L^p}.
 \end{align} 
In the following subsection, we estimate the part of the series on the right-hand side of~\eqref{E:series} in which $j> 0$, and in the third subsection we consider the case when $j \le 0$.

\subsection{The case $j>0$}\label{proof3}

The estimate derived in this subsection will be based on a combination of Theorem~\ref{T:HormanderB} and of the following proposition. Its variant assuming $p_1<\infty$ and $p_2 <\infty$ was proved in \cite[Lemma 11]{Gr_He_Ho2018}, and the extension to the case $p_1=\infty$ or $p_2=\infty$ was obtained in~\cite[(3.3)]{He_Park2022}.

\begin{customproposition}{B}[\cite{Gr_He_Ho2018, He_Park2022}]\label{pr:Linfty_bound}
Let $1<p_1,p_2\leq \infty$ and $1\leq p<\infty$ satisfy $1/p=1/p_1+1/p_2$. Suppose that $j> 0$.
Then there exists a constant $\delta_0>0$, depending on $p_1,p_2$, such that
\begin{equation*}
\Vert T^{j}\Vert_{L^{p_1}\times L^{p_2}\to L^p}\lesssim_{n,\delta_0} 2^{-\delta_0 j}\Vert \Omega\Vert_{L^\infty(\bbs^{2n-1})}.
\end{equation*}
\end{customproposition}%

We next bound $\sum_{j=1}^\infty \Vert T^{j}\Vert_{L^{p_1}\times L^{p_2}\to L^p}$.
Without loss of generality, we may assume that $\|\Omega\|_{L(\log L)^A(\bbs^{2n-1})} =1$.  
We decompose the sphere $\bbs^{2n-1}$ as 
$$\bbs^{2n-1}=\dot{\bigcup_{\mu\in\bbn_0}}D^{\mu}$$
where
\begin{equation*}
		D^\mu:=\begin{cases}
		\big\{\theta\in\S^{2n-1}\,:\, |\Omega(\theta)|\leq 1 \big\}  & \text{if }\mu=0,\\
		 \big\{\theta\in\S^{2n-1}\,:\, 2^{\mu-1} < |\Omega(\theta)|\le 2^\mu\big\}  &\text{if }\mu\geq 1,
		 \end{cases}
	\end{equation*}
and write
$$\Omega(\theta)=\Omega(\theta)-\int_{\bbs^{2n-1}}\Omega(\eta)     \;d\nu(\eta) = \sum_{\mu=0}^{\infty}\bigg( \Omega(\theta)\chi_{D^{\mu}}(\theta)-\int_{D^{\mu}} \Omega(\eta) \;d\nu(\eta)\bigg)=:\sum_{\mu=0}^{\infty}\Omega^{\mu}(\theta).$$
Then each $\Omega^{\mu}$ preserves the vanishing moment condition 
$$\int_{\bbs^{2n-1}}\Omega^{\mu}(\theta)\;d\nu(\theta)=0.$$
Moreover, it also has the norm estimates
\begin{equation*}
\Vert \Omega^{\mu}\Vert_{L^{\infty}(\bbs^{2n-1})}\le 2^{\mu+1}
\end{equation*}	
	and 
\begin{align*}
\|\Omega^\mu\|_{L^1(\bbs^{2n-1})}& \le2\int_{D_\mu}|\Omega(\theta)| \;d\nu(\theta).
\end{align*}

Let $T^{j,\mu}$ and $K^{j,\mu}_0$ be the counterparts of $T^j$ and $K^j_0$ when we replace $\Omega$ by $\Omega^{\mu}$.
 Then Proposition \ref{pr:Linfty_bound} yields that
\begin{equation}\label{lowermuest}
\big\Vert T^{j,\mu}\big\Vert_{L^{p_1}\times L^{p_2}\to L^p}\lesssim 2^{-\delta_0 j} 2^{\mu}.
\end{equation}
Furthermore, we see that
\begin{equation}\label{E:FSupp}\supp\big(\wh{K^{j,\mu}_0}\big)\subset \big\{\xxxi\in \bbr^{2n}: 2^{j-1}\le |\xxxi|\le 2^{j+1} \big\}\end{equation}
and the operator $T^{j,\mu}$ is of the form~\eqref{E:b} with $K$ replaced by $(K^{j,\mu}_0)_{-j}$. 
Thus, with reference to Remark \ref{remark1},
 Theorem~\ref{T:HormanderB} implies
\begin{align}\label{tmujpppest}
\big\Vert T^{j,\mu}\big\Vert_{L^{p_1}\times L^{p_2}\to L^p}&\lesssim 2^{-2jn}\int_{\R^{2n}} \big|K^{j,\mu}_0(2^{-j}\yyy)\big|(\log(e+|\yyy|))^{A-1} \,d\yyy\nonumber\\
&\lesssim j^{A-1}\int_{\R^{2n}} \big|K^{j,\mu}_0(\yyy)\big|(\log(e+|\yyy|))^{A-1}\,d\yyy
\end{align}
for any real number $A$ satisfying $A \geq \max\{\frac{1}{p_1},\frac{1}{p_2},\frac{1}{p'}\}+1$.
Let $B(0,R)$ denote the ball of radius $R$, centered at the origin in $\bbr^{2n}$.
We use the superposition estimate
$$\big| \Gamma(\yyy)\big|\lesssim\sum_{r=1}^{\infty}2^{-4nr}\chi_{B(0,2^r)}(\yyy),$$
which yields
$$\big| \Gamma_j(\yyy)\big|\lesssim\sum_{r=1}^{\infty}2^{2nj-4nr}\chi_{B(0,2^{r-j})}(\yyy).$$
Moreover, using \eqref{E:FSupp}, we can see that $$|K^{j,\mu}_0(\xxx)|\lesssim \sum_{l=j-2}^{j+2} |K^{j,\mu}_0(\xxx)|\ast |\Gamma_l(\xxx)|.$$ The integral on the right-hand side of \eqref{tmujpppest} can therefore be controlled by
\begin{align*}
& \sum_{l=j-2}^{j+2} \sum_{r=1}^{\infty}2^{2nl-4nr}\int_{\bbr^{2n}}|K^{j,\mu}_0|\ast \chi_{B(0,2^{r-l})}(\yyy)(\log{(e+|\yyy|))^{A-1}} \;d\yyy\\
&\lesssim \sum_{r=1}^{\infty}r^{A-1}2^{2nj-4nr}\big\Vert |K^{j,\mu}_0|\ast \chi_{B(0,2^{r-j})} \big\Vert_{L^1(\bbr^{2n})}\\
&\lesssim \Vert K^{j,\mu}_0\Vert_{L^1(\bbr^{2n})}\sum_{r=1}^{\infty}r^{A-1}2^{-2nr}\lesssim \Vert K^{j,\mu}_0\Vert_{L^1(\bbr^{2n})}\lesssim \Vert \Omega^{\mu}\Vert_{L^1(\mathbb{S}^{2n-1})}.
\end{align*}
This finally yields
\begin{equation}\label{uppermuest}
\big\Vert T^{j,\mu}\big\Vert_{L^{p_1}\times L^{p_2}\to L^p}\lesssim j^{A-1} \int_{D_\mu}|\Omega(\theta)|\;d\nu(\theta).
\end{equation}

Now, choosing $0<\epsilon<\delta_0$, we decompose
\begin{equation*}
\sum_{j=1}^{\infty}\big\Vert T^j\big\Vert_{L^{p_1}\times L^{p_2}\to L^p}\le \sum_{j=1}^{\infty}\sum_{0\le \mu<j\epsilon}\big\Vert T^{j,\mu}\big\Vert_{L^{p_1}\times L^{p_2}\to L^p}+\sum_{j=1}^{\infty}\sum_{\mu\ge j\epsilon}\big\Vert T^{j,\mu}\big\Vert_{L^{p_1}\times L^{p_2}\to L^p}.
\end{equation*}
By~\eqref{lowermuest}, the first sum is controlled by
\begin{align*}
\sum_{j=1}^{\infty}2^{-\delta_0 j}\sum_{0\le \mu<j\epsilon}2^{\mu}\lesssim \sum_{j=1}^{\infty}2^{-j(\delta_0-\epsilon)}\lesssim 1.
\end{align*}
For the remaining term, we apply \eqref{uppermuest} to obtain
\begin{align*}
\sum_{j=1}^\infty \sum_{\mu\ge j\epsilon}\big\Vert T^{j,\mu}\big\Vert_{L^{p_1}\times L^{p_2}\to L^p}
&\lesssim  \sum_{j=1}^\infty j^{A-1}\sum_{\mu\ge j\epsilon} \int_{D^\mu}|\Omega(\theta)|\;d\nu(\theta)\\
&=\sum_{\mu=1}^\infty \int_{D^\mu}|\Omega(\theta)|\;d\nu(\theta) \sum_{j \leq \mu/\epsilon} j^{A-1}\\
&\lesssim \sum_{\mu=1}^\infty \mu^A \int_{D^{\mu}}|\Omega(\theta)|\;d\nu(\theta)\\
&\lesssim \int_{\S^{2n-1}} |\Omega(\theta)| (\log(e+|\Omega(\theta)|))^{A}\;d\nu(\theta)\lesssim 1,
\end{align*}
as desired.

\subsection{The case $j\le 0$}

In this subsection, we apply a bilinear version of Mikhlin's multiplier theorem, which was established by Coifman and Meyer \cite{Co_Me1978} under the assumption $1<p<\infty$ and later extended to the case $p \leq 1$ by Kenig and Stein \cite{Ke_St1999} and Grafakos and Torres \cite{Gr_To2002}.
\begin{customlemma}{C}[\cite{Co_Me1978, Gr_To2002,  Ke_St1999}]\label{CGKlemma}
Let $1< p_1,p_2\le \infty$ and $1/2< p< \infty$ with $1/p_1+1/p_2=1/p$.
Suppose that a bounded function $m$ on $\bbr^{2n}$ satisfies
\begin{equation*}
\big| \partial_{\xxxi}^{\vec{\bold \alpha}}m(\xxxi)\big|\lesssim_{\vec{\alpha}}\mathfrak{A}|\xxxi|^{-|\vec{\alpha}|}
\end{equation*} for all multi-indices $\vec{\alpha}\in \bbn_0^{2n}$ with $|\vec{\alpha}|\le 2n$.
Then the corresponding bilinear Fourier multiplier operator $\mathcal{B}$ defined by~\eqref{E:bilinear-fm} is bounded from $L^{p_1}(\R^n)\times L^{p_2}(\R^n)$ to $L^p(\R^n)$ and 
$$\big\Vert \mathcal B\big\Vert_{L^{p_1}\times L^{p_2}\to L^p}\lesssim \mathfrak{A}.$$
\end{customlemma} 

We note that the operator $T^j$ defined in \eqref{tjdefinition} coincides with the bilinear Fourier multiplier operator associated with the symbol $m=\wh{K^j}$, and it thus suffices to prove 
\begin{equation}\label{kjmultipliercondition}
\big| \partial_{\xxxi}^{\vec{\alpha}}\wh{K^j}(\xxxi)\big|\lesssim_{\vec{\alpha}}2^j \Vert \Omega\Vert_{L^1(\bbs^{2n-1})} |\xxxi|^{-|\vec{\alpha}|},\qq j\le 0
\end{equation}
for multi-indices $\vec{\alpha}$ with $|\vec{\alpha}|\le 2n$.
Once inequality~\eqref{kjmultipliercondition} is established, an application of Lemma~\ref{CGKlemma} will yield
\[
\|T^j\|_{L^{p_1}\times L^{p_2}\to L^p} \lesssim 2^j \Vert \Omega\Vert_{L^1(\bbs^{2n-1})}
\lesssim 2^j \Vert \Omega\Vert_{L(\log L)^A(\bbs^{2n-1})}, \quad j \leq 0,
\]
which is summable in $j$, as desired. 

The estimate \eqref{kjmultipliercondition} can be verified by mimicking the argument in \cite[Corollary 4.1]{Du_Ru1986}, which deals with the case when $\Omega$ belongs to $L^q(\bbs^{mn-1})$ for $q>1$. We include the proof of \eqref{kjmultipliercondition} for the sake of completeness.

We first claim that
\begin{equation}\label{dek0est}
\big| \partial_{\xxxi}^{\vec{\alpha}}\wh{K_0}(\xxxi)\big|\lesssim_{\vec{\alpha}}\Vert \Omega\Vert_{L^1(\bbs^{2n-1})}\begin{cases}
|\xxxi| & \text{ if }~\vec{\alpha}=\vec{0},\\
1 & \text{ if }~\vec{\alpha}\not= \vec{0}.
\end{cases}
\end{equation}
Indeed, by applying the vanishing moment condition \eqref{E:vanishingcon}, we write
\begin{align*}
\wh{K_0}(\xxxi)&=\int_{\bbs^{2n-1}}\Omega(\theta)\bigg( \int_{1/2}^{2}e^{-2\pi i r\langle \theta,\xxxi \rangle} \wh{\Gamma}(r\theta)  \; \frac{dr}{r}\bigg) \;d\nu(\theta)\\
&=\int_{\bbs^{2n-1}}\Omega(\theta)\bigg( \int_{1/2}^{2}\big( e^{-2\pi i r\langle \theta,\xxxi \rangle} -1\big) \wh{\Gamma}(r\theta)  \; \frac{dr}{r}\bigg) \;d\nu(\theta),
\end{align*}
and thus
$$\big| \wh{K_0}(\xxxi)\big|\le \int_{\bbs^{2n-1}} \big| \Omega(\theta)\big| \big| \langle \theta,\xxxi\rangle\big|\bigg(\int_{1/2}^{2}\big|\wh{\Gamma}(r\theta) \big|\frac{dr}{r}  \bigg)  \; d\nu(\theta)\lesssim \Vert \Omega\Vert_{L^1(\bbs^{2n-1})}|\xxxi|,$$
as desired. Moreover, for multi-indices $\vec{\alpha}\not=\vec{0}$,
\begin{align*}
\big| \partial_{\xxxi}^{\vec{\alpha}}\wh{K_0}(\xxxi)\big|&=\bigg| \int_{\R^{2n}}  (-2\pi i \xxx\big)^{\vec{\alpha}}K_0(\xxx)e^{2\pi i\langle \xxx,\xxxi\rangle}   \; d\xxx \bigg|\\
&\le \big( 2\pi\big)^{|\vec{\alpha}|} \int_{\bbs^{2n-1}}\big| \Omega(\theta)\big| \bigg( \int_{1/2}^{2}r^{|\vec{\alpha}|-1} \big| \wh{\Gamma}(r\theta)\big|\; dr\bigg) \; d\nu(\theta)\\
&\lesssim \Vert \Omega\Vert_{L^1(\bbs^{2n-1})},
\end{align*}
which completes the proof of \eqref{dek0est}.

Therefore, in view of \eqref{whkjxiex}, we have
\begin{align*}
\big| \wh{K^j}(\xxxi)\big|&\le \sum_{k\in\bbz} \big| \wh{K_0}(2^{-k}\xxxi)\big|\chi_{|\xxxi|\sim 2^{j+k}}(\xxxi)\\
&\lesssim \Vert \Omega\Vert_{L^1(\bbs^{2n-1})}\sum_{k\in\bbz}2^{-k} |\xxxi|\chi_{|\xxxi|\sim 2^{j+k}}(\xxxi)\lesssim 2^j \Vert \Omega\Vert_{L^1(\bbs^{2n-1})},
\end{align*}
which proves \eqref{kjmultipliercondition} for $\vec{\alpha}=\vec{0}$.
If $\vec{\alpha}\not= \vec{0}$,
then
\begin{align*}
\big| \partial_{\xxxi}^{\vec{\alpha}}\wh{K^j}(\xxxi)\big|&\lesssim \sum_{k\in\bbz}\sum_{\vec{\beta}+\vec{\gamma}=\vec{\alpha}}2^{-k|\vec{\beta}|}\big| \partial_{\xxxi}^{\vec{\beta}}\wh{K_0}(2^{-k}\xxxi)\big|2^{(-k-j)|\vec{\gamma}|}\chi_{|\xxxi|\sim 2^{j+k}}(\xxxi)\\
&\lesssim \Vert \Omega\Vert_{L^1(\bbs^{2n-1})}\sum_{k\in\bbz}\big| 2^{-k}\xxxi\big|2^{(-k-j)|\vec{\alpha}|}\chi_{|\xxxi|\sim 2^{j+k}}(\xxxi)\\
&\qq+\Vert \Omega\Vert_{L^1(\bbs^{2n-1})}\sum_{k\in\bbz}\sum_{\vec{\beta}+\vec{\gamma}=\vec{\alpha},\vec{\beta}\not=\vec{0}}2^{-k|\vec{\beta}|}2^{(-k-j)|\vec{\gamma}|}\chi_{|\xxxi|\sim 2^{j+k}}(\xxxi)\\
&\lesssim 2^j \Vert \Omega\Vert_{L^1(\bbs^{2n-1})} |\xxxi|^{-|\vec{\alpha}|}\sum_{k\in\bbz}\chi_{|\xxxi|\sim 2^{j+k}}(\xxxi)\\
&\sim 2^j \Vert \Omega\Vert_{L^1(\bbs^{2n-1})} |\xxxi|^{-|\vec{\alpha}|},
\end{align*}
where the penultimate estimate holds because $j\le 0$. This completes the proof of \eqref{kjmultipliercondition}.

\section*{Acknowledgments}

We thank the anonymous referees for their comments which helped to improve the exposition.
We also thank Tengfei Bai for pointing out to us an oversight in the proof of Theorem 3.


\begin{thebibliography}{99}

\bibitem{BS}
A. Baernstein and E. T. Sawyer, \emph{Embedding and multiplier theorems for $H^p(R^n)$}, Mem. Amer. Math. Soc. \textbf{53} (1985), no. 318.

\bibitem{CT}  A. P.  Calder\'on and A. Torchinsky,
  {\it Parabolic maximal functions associated with a distribution,  II,}
  Adv. Math.  \textbf{24} (1977)  101--171.

\bibitem{Ca_Zy1956}
A. P. Calder\'on and A. Zygmund, \emph{On singular integrals}, Amer. J. Math. \textbf{78} (1956), 289-309.


\bibitem{Car}
A.  Carbery, \emph{Variants of the Calderón-Zygmund theory for $L^p$--spaces}, Rev. Mat. Iberoam. \textbf{2} (1986), no. 4, 381--396.

\bibitem{Ch_Ru1988}
M. Christ and J.-L. Rubio de Francia,  \emph{Weak type $(1,1)$ bounds for rough operators II}, Invent. Math. \textbf{93} (1988) 225-237.



\bibitem{Co_Me1978}
R. R. Coifman and Y. Meyer,  \emph{Au del\`a des op\'erateurs pseudo--diff\'erentiels}, Ast\'erisque \textbf{57}
(1978) 1--185.



\bibitem{Do_Sl_submitted}
G. Dosidis and L. Slav\'ikov\'a,  \emph{Multilinear singular integrals with homogeneous kernels near $L^1$}, Math. Ann. \textbf{389} (2024), no. 3, 2259--2271.

\bibitem{Du_Ru1986}
J. Duoandikoetxea and J. L. Rubio de Francia,  \emph{Maximal and singular integral operators via Fourier transform estimates}, Invent. Math. \textbf{84} (1986) 541--561.
	

\bibitem{Fe_St1971}
C. Fefferman and E. M. Stein,  \emph{Some maximal inequalities}, Amer. J. Math. \textbf{93} (1971) 107--115.


\bibitem{Fr_Ja1990}
M. Frazier and B. Jawerth,  \emph{A discrete transform and decomposition of distribution spaces}, J. Funct. Anal. \textbf{93}
(1990) 34--170.

\bibitem{MFA}
L. Grafakos,
\newblock {\em Modern Fourier Analysis.,}  3rd Ed., GTM {\bf 250}, 
\newblock Springer, 2014.

\bibitem{Gr_He_Ho2018}
L. Grafakos, D. He, and P. Honz\'ik, \emph{Rough bilinear singular integrals}, Adv. Math. {\bf 326} (2018) 54--78.

\bibitem{Gr_He_Ho_Park_submitted}
L. Grafakos, D. He, P. Honz\'ik, and B. Park,  \emph{Multilinear rough singular integral operators}, J. Lond. Math. Soc. (2) \textbf{109} (2024), no. 2, Paper No. e12867, 35 pp.

\bibitem{GHS20}
L.~Grafakos, D.~He, and L.~Slav{\'i}kov{\'a},
\newblock\emph{$L^2 \times L^2 \rightarrow L^1$ boundedness criteria},
\newblock Math. Ann. {\bf 376} (2020), no. 1--2, 431--455.

\bibitem{Gr_To2002}
L. Grafakos and R. Torres,  \emph{Multilinear Calder\'on-Zygmund theory}, Adv. Math. \textbf{165} (2002) 124--164.

\bibitem{Haar}
A. Haar, \emph{The Interplay of Shifted Square and Maximal Function Estimates in the Context of Multilinear Fourier Multipliers}, preprint, arXiv:2512.01599.

\bibitem{He_Park2022}
D. He and B. Park,  \emph{Improved estimates for bilinear rough singular integrals}, Math. Ann. \textbf{386} (2023) 1951-1978.

\bibitem{Ho1988}
S. Hofmann, \emph{Weak type $(1,1)$ boundedness of singular integrals with nonsmooth kernels}, Proc. Amer. Math. Soc. \textbf{103} (1988) 260-264.


\bibitem{Hoe}  L. H\"ormander,
\emph{Estimates for translation invariant operators
in $L^p$ spaces,} Acta Math. {\bf 104} (1960) 93--140.


\bibitem{Ke_St1999}
C. Kenig and E. M. Stein,  \emph{Multilinear estimates and fractional integration}, Math. Res. Lett.  \textbf{6} (1999) 1--15.


\bibitem{Mu2014} 
C. Muscalu, \emph{Calder\'on Commutators and the Cauchy integral on Lipschitz curves revisited: I. First commutator and generalizations}, Rev. Mat. Iberoam. \textbf{30} (2014) 1413--1437.


\bibitem{Park2019}
B. Park,  \emph{Some maximal inequalities on Triebel--Lizorkin spaces for $p=\infty$},  Math. Nachr. \textbf{292} (2019) 1137--1150.

\bibitem{Park2019_2}
B. Park,  \emph{Fourier multiplier theorems for Triebel--Lizorkin spaces}, Math. Z. \textbf{293}
(2019) 221--258.

\bibitem{Park_IUMJ}
B. Park,  \emph{Equivalence of (quasi--)norms on a vector--valued function space and its applications to multilinear operators}, Indiana Univ. Math. J. \textbf{70} (2021) 1677--1716.



\bibitem{Pe1975}
J. Peetre,  \emph{On spaces of Triebel-Lizorkin type}, Ark. Mat. \textbf{13}
(1975) 123-130.	

\bibitem{Seeger88}
A. Seeger, \emph{Some inequalities for singular convolution operators in $L^p$--spaces},
Trans. Amer. Math. Soc. \textbf{308} (1988), no. 1, 259--272.
	
\bibitem{Se1990}
A. Seeger,  \emph{Remarks on singular convolution operators}, Studia Math. \textbf{97}
(1990), no. 2, 91--114.	

\bibitem{Se1996}
A. Seeger, \emph{Singular integral operators with rough convolution kernels}, J. Amer. Math. Soc. \textbf{9} (1996) 95-105.

	
\bibitem{St1993}
E. M. Stein,  \emph{Harmonic Analysis, Real Variable Methods, Orthogonality, and Oscillatory Integrals}, Princeton University Press, 1993.

\bibitem{Tomita}
N. Tomita, \emph{On the Hörmander multiplier theorem and modulation spaces}, Appl. Comput. Harmon. Anal. \textbf{26} (2009), no. 3, 408--415.

\bibitem{Tr}
H. Triebel, \emph{Theory of Function Spaces}, Birkhauser, Basel-Boston-Stuttgart
(1983).


	
	
\end{thebibliography}
\end{document}